\documentclass[11pt]{article}

\usepackage{booktabs}
\usepackage[margin=15pt,font=small,labelfont={bf,sf},justification=justified]{caption}
\usepackage[superscript,biblabel]{cite}
\usepackage{color}
\usepackage{float}
\usepackage{graphicx}
\usepackage{lscape}
\usepackage[bbgreekl]{mathbbol}
\usepackage{mathtools}
\usepackage{multibib}
\usepackage{multirow}
\usepackage{tabularx}
\usepackage{titleref}
\usepackage{url}

\usepackage{tikz}
\usepackage{amsmath}
\usepackage{subcaption}

\usetikzlibrary{fadings}
\usetikzlibrary{patterns}
\usetikzlibrary{shadows.blur}
\usetikzlibrary{shapes}

\newif\ifbembo
\newif\ifcharter
\newif\iferewhon
\newif\iflibertine
\newif\iflibertinealt
\newif\ifpalantino
\newif\iftimesnewroman

\bembofalse
\chartertrue
\erewhonfalse
\libertinefalse
\libertinealtfalse
\palantinofalse
\timesnewromanfalse
\ifbembo
  \usepackage[p,osf]{ETbb} % osf in text, tabular lining figures in math
  \usepackage[scaled=.95,type1]{cabin} % sans serif in style of Gill Sans
  \usepackage[varqu,varl]{zi4}% inconsolata typewriter
  \usepackage[T1]{fontenc}
  \usepackage[libertine,vvarbb]{newtxmath}
  \usepackage[cal=boondoxo,bb=boondox]{mathalfa}
\fi

\ifcharter
  \usepackage[scaled=.98,p]{XCharter}
  \usepackage[scaled=1.04,varqu,varl]{inconsolata}
  \usepackage[type1]{cabin}
  \usepackage[uprightscript,xcharter,vvarbb,scaled=1.05]{newtxmath}
  \usepackage[cal=euler,scr=boondoxo,bb=boondox]{mathalpha}  % prefer bb=bboldxLight when available
\fi

\iferewhon
  \usepackage[p,osf,scaled=.98,space]{erewhon} % scaling by .98 not really necessary
  \usepackage[varqu,varl]{inconsolata} % typewriter
  \usepackage[type1,scaled=.95]{cabin} % sans serif like Gill Sans
  \usepackage[utopia,vvarbb]{newtxmath}
\fi

\iflibertine
  \usepackage[sb]{libertinus}
  \usepackage[T1]{fontenc}
  \usepackage{textcomp}
  \usepackage[varqu,varl]{zi4}% inconsolata for mono, not LibertineMono
  \usepackage{libertinust1math} % slanted integrals, by default
  \usepackage[scr=boondoxo,bb=boondox]{mathalpha} %Omit bb=boondox for default libertinus bb
\fi

\iflibertinealt
  \usepackage[sb]{libertinus}
  \usepackage[T1]{fontenc}
  \usepackage{textcomp}
  \usepackage{libertinust1math} % slanted integrals, by default
  \usepackage[scr=boondoxo,bb=boondox]{mathalpha} %Omit bb=boondox for default libertinus bb
\fi

\ifpalantino
  \usepackage[largesc]{newpxtext}
  \usepackage[vvarbb]{newpxmath}
\fi

\iftimesnewroman
  \usepackage{newtxtext} %T1 is default encoding
  \usepackage[scaled=0.95]{inconsolata} % typewriter
  \usepackage[vvarbb]{newtxmath} % vvarbb gives STIX Bbb
\fi

\usepackage[T1]{fontenc}
\usepackage[final,protrusion=true,expansion=true]{microtype}

\usepackage{titling}
\usepackage{authblk} % note: titling needs to be loaded before authblk
\pretitle{\begin{center}\bfseries\sffamily\LARGE}
\posttitle{\end{center}}

\predate{\begin{center}\normalsize}
\postdate{\end{center}}

\usepackage{abstract}
\usepackage{sectsty} % note: the main alternative package, titlesec, does not work with titleref
\allsectionsfont{\sffamily}

\usepackage[left=1in,right=1in,top=1in,bottom=1in,includefoot,heightrounded]{geometry}

\newcites{supp}{SUPPLEMENTAL REFERENCES}

\makeatletter
\patchcmd{\LS@rot}{90}{-90}{}{}
\patchcmd{\endlandscape}{90}{-90}{}{}
\makeatother
\newcommand{\Rey}{\mathrm{Re}}

\newcommand{\bF}{\mathbf{F} }

\newcommand{\bP}{\boldsymbol P}

\newcommand{\bphi}{\boldsymbol{\phi}}
\newcommand{\bpsi}{\boldsymbol{\psi}}

\newcommand{\bn}{\mathbf n}

\newcommand{\Tau}{\mathrm{T}}

\newcommand{\bu}{\boldsymbol u}
\newcommand{\bU}{\mathbf U}

\newcommand{\norm}[1]{\left\lVert #1 \right\lVert}

\DeclareGraphicsExtensions{.pdf,.eps,.ps,.eps.gz,.ps.gz,.eps.Y}

\newcommand{\bchi}{\boldsymbol{\chi}}

\newcommand{\compositenorm}[1]{{\left\vert\kern-0.25ex\left\vert\kern-0.25ex\left\vert #1
    \right\vert\kern-0.25ex\right\vert\kern-0.25ex\right\vert}}

\newcommand{\jumpinverse}[1]{\llbracket#1\rrbracket}

\newcommand*{\tran}{^{\mkern-1.5mu\mathsf{T}}}

\begin{document}

\title{A Pressure-Robust Immersed Interface Method for Discrete Surfaces}

\author[1]{Michael J. Facci$^*$}
\author[1]{Qi Sun\footnote{These two authors contributed equally to this work.}}
\author[1,2--5]{Boyce E.~Griffith\thanks{Corresponding author: boyceg@email.unc.edu}}
\affil[1]{Department of Mathematics, University of North Carolina, Chapel Hill, NC, USA}
\affil[2]{Department of Biomedical Engineering, University of North Carolina, Chapel Hill, NC, USA}
\affil[3]{Carolina Center for Interdisciplinary Applied Mathematics, University of North Carolina, Chapel Hill, NC, USA}
\affil[4]{Computational Medicine Program, University of North Carolina School of Medicine, Chapel Hill, NC, USA}
\affil[5]{McAllister Heart Institute, University of North Carolina School of Medicine, Chapel Hill, NC, USA\vspace{\baselineskip}}
\maketitle

\begin{abstract}

The immersed interface method (IIM) for fluid-structure interaction imposes discontinuities in the fluid stress along immersed boundaries that are generated by forces concentrated along those boundaries. For a viscous incompressible fluid, imposing these discontinuities requires decomposing the boundary force into its normal and tangential components, which respectively determine jump conditions for the pressure and velocity gradient. In previous work, we developed an IIM for geometries described by $C^0$ triangulated surfaces, with a focus on piecewise linear surface representations. In this setting, the normal and tangent vectors of the discrete surface are constant on each element, and that method uses those piecewise constant vectors to determine the normal and tangential force components and, ultimately, the jump condition.
We demonstrated that this is substantially more accurate than immersed boundary methods that use regularized delta functions at corresponding grid resolutions for situations in which shear stresses dominate, as in external flow applications.
However, this IIM formulation struggles to accurately capture pressure loads. Here, we identify that the primary cause of this limitation is the discontinuous surface normal inherent in $C^0$ triangulated surfaces.
To address this, we propose a procedure that uses approximations of the surface normals that more accurately account for the curvature and avoid discontinuities in the reconstructed normal and tangent vectors.
%Two formulations for the discrete projection are considered: a nodal quadrature method, corresponding to linear interpolation of nodal normal vectors, and a consistent formulation based on the finite element basis functions.
In this paper, we investigate two ways to construct the continuous surface normal approximation. First, we construct a continuous approximation by performing a standard $L^2$ projection of the discontinuous surface normal field into a continuous finite element space. Second, we construct vertex normal vectors using inverse centroid-distance weighting and apply linear interpolation to define a continuous normal field.
Numerical experiments demonstrate that the use of jump conditions computed with reconstructed continuous normal vector fields reduce leakage by up to six orders of magnitude across a range of pressure loads. This work offers a major improvement to the volume conservation properties of our IIM formulation, thereby facilitating its application to models involving both large sheer stresses and pressure loads.
%Normal components of the force generate pressure discontinuities, and tangential components generate jumps in the viscous shear stress. Evaluating these discontinuities thereby requires descriptions of the normal and tangent vectors along the interface. In earlier work, we developed an IIM for geometries described by $C^0$ triangulated surfaces. Normal and tangent vectors were directly evaluated from the piecewise polynomial elements that comprise the surface triangulation. We demonstrated that this method can be highly effective for external and internal flow applications that do not include large pressure differences across the interface. This methodology can provide poor volume conservation in applications involving large pressure loads, however. Herein, we identify that surface normal vector construction, which is discontinuous since the geometrical description is only $C^0$, as the primary cause of this poor volume conservation. 
%To improve the volume conservation of the methodology, we introduce two straightforward approaches to construct the projected vector field. The first computes the $L^2$ projection of the geometrical normal vectors, which are discontinuous at element boundaries when using standard $C^0$ Lagrangian elements, onto the continuous Lagrange basis functions on the interface. The second approach projects the discontinuous vectors onto the same continuous space by employing a lumped mass matrix, a procedure consistent with Phong shading, a computer graphics technique that linearly interpolates vertex normals across discrete element faces. 

\end{abstract}


\section{Introduction}
Accurately and efficiently resolving immersed boundaries and their interactions with fluids remains a central challenge in numerical simulations of fluid–structure interaction (FSI). Traditional body-fitted grid methods~\cite{doi:https://doi.org/10.1002/0470091355.ecm009,SOULI2000659,hou2012numerical}, whether structured or unstructured, can offer high-order accuracy for capturing boundary dynamics. However, these methods can become computationally expensive and difficult to implement when dealing with large deformations or complex moving geometries because they require  continuous grid regeneration. 

To overcome these limitations, a variety of methods that avoid body-fitted grids have been developed, which reduce the need of dynamic grid adaptation and thereby improve computational efficiency and simplify implementation.  Among them are Peskin's immersed boundary method,~\cite{PESKIN1972252,Peskin_2002,ADJERID2015170}  the extended finite element method~\cite{CLAUS2019185,He_Song2019,Massing_Larson2014,Wang_Chen2019}, the immersed interface method (IIM)~\cite{li2003overview, XU2006454,lee2003immersed,tan2009}, and other related approaches~\cite{TSENG2003593,OlshanskiiQiQuaini2021,OlshanskiiQiQuaini2022}.
Classical IB methods have substantial pressure leaks, although these can be mitigated using composite B-spline regularized delta functions\cite{gruninger2026composite} and other related methods\cite{bao2017immersed,CHEN2024113048}.

Classical IIMs~\cite{li2003overview, XU2006454,lee2003immersed,tan2009} were developed using at smooth (at least $C^2$) representations of the interface geometry. This is because the key idea of IIMs is to incorporate the known jump conditions of the analytic solution and its derivatives into the background fluid discretization schemes near the discontinuities caused by the singular force of the immersed interface. For example, one can derive the principal jump conditions~\cite{XU2006454} for pressure balance from the normal components of the interfacial force, and the tangential components produce discontinuities in the viscous shear stress. However, specializing these methods to $C^2$ geometry limits the applicability of the IIM to real-world problems in engineering and applied science, because typical mesh generation for complex geometries often provides only a $C^0$ representation and may contain sharp corners or edges even if the geometry of interest is itself smooth. 
To address this limitation, we developed an IIM framework~\cite{kolahdouz_immersed_2020,FACCI2025114119,sun2025} capable of handling interfaces that are described only in terms of a surface triangulation, without requiring the underlying analytical surface description. This formulation avoids the need for globally smooth geometry and supports arbitrary polyhedral meshes, requiring only a discrete representation of the surface. Prior work demonstrated that this approach is effective for external flow problems, offering high spatial accuracy and geometric flexibility~\cite{KOLAHDOUZ2023112174}.
Building on this formulation, we constructed a coupled FSI framework~\cite{KOLAHDOUZ2023112174,KOLAHDOUZ2021110442} that accurately resolves interactions between fluids and both rigid and deformable structures. This method has been successfully applied to problems involving external flows around flexible membranes and oscillating plates. Biomedical applications have also been presented, including simulations of the dynamics of a bileaflet mechanical heart valve in a pulse duplicator system and the transport of blood clots in a patient-specific anatomical model of the inferior vena cava. These simulations show strong agreement with analytical benchmarks and experimental data.

Applying this approach to internal flows characterized by significant pressure loading, such as those encountered in pulsatile blood flow, exposes a critical numerical issue: the method struggles to enforce the no-penetration condition along the immersed boundary.
Spurious fluxes across the interface arise in models involving large pressure loads. In that method, the jump conditions are computed using the discrete surface normal, which deviate from the true jump condition because the discrete $C^0$ geometry does not capture the curvature nor the continuous normal vector field of the true geometry.  These errors in the discrete jump conditions lead to non-physical flows through the boundary with magnitudes that are determined by the pressure jump. Although these spurious fluxes can be eliminated through grid refinement, situations that involve very large pressure loads require prohibitively high spatial resolution.

\iffalse
This deficiency is most pronounced near corners or edges of discrete surfaces, where discontinuities in the discrete surface normal vectors are significant. These discontinuities introduce aliasing errors and result in localized fluid leakage. Although the leakage diminishes with grid refinement—particularly when the analytical geometry is $C^1$ smooth—it does so slowly. As a result, the computational cost required to reduce the leakage can quickly accumulate. If left unaddressed, the cumulative effect may significantly distort simulation outcomes, especially under strong pressure gradients.
\fi
The key contribution of this study is that it introduces a new formulation for the projected interface condition based on a smoothed reconstruction of the surface normal vector field. In particular, we focus on piecewise linear surface representations.
Smoothed reconstructions of a piecewise linear surface's curvature and normal vectors are often studied in the field of computer graphics. One such example in Chen et al.~\cite{CHEN2004447} demonstrates that smoothing the discrete normal vector using inverse centroid-weights substantially decreases error compared to area-weighted normals and flat normals in the computed surface curvature and vertex normal vectors. Following a similar approach, we present two smoothing methods to improve the accuracy of the reconstructed normal vectors and surface curvature, and investigate their effects on numerical leakage.
In one approach, we project the discontinuous geometrical surface normals onto a continuous finite element space.
In another approach, we use an inverse centroid-distance weighting scheme~\cite{CHEN2004447} to compute vertex normals, which are then linearly interpolated~\cite{phong} to provide a continuous representation across the interface.
We demonstrate that by computing our jump conditions using smoothed normal vector reconstructions, we achieve substantial improvement in enforcing the non-penetration condition while retaining a piecewise linear surface representation. In our tests, we achieve up to six orders of magnitude reduction in spurious flux across a range of pressure-driven internal flow problems, without otherwise degrading spatial accuracy or robustness. These results demonstrate that the proposed formulation provides both flexibility and improved conservation properties, making it particularly suitable for large-scale, three-dimensional simulations involving internal flows under significant pressure loads.

\section{Continuous Equations of Motion}

This section outlines the equations of motion for a viscous incompressible fluid with interfacial forces concentrated along an immersed boundary. This boundary both applies force to the fluid and moves with the local fluid velocity. The fluid domain, $\Omega $, is partitioned
into two subregions, $\Omega_t^+$ (exterior) and $\Omega_t^-$ (interior), by an
interface $\Gamma_t = \overline{\Omega_t^+} \bigcap \overline{\Omega_t^-}$.  We use the initial coordinates $\mathbf{X}$ of the boundary as Lagrangian coordinates, so that $\mathbf{X} \in \Gamma_0$, with corresponding current coordinates $\boldsymbol{\chi}(\mathbf{X},t) \in \Gamma_t$. 
The equations of motion are
\begin{align}
  \rho \frac{\mathrm{D} \mathbf{u}(\mathbf{x},t)}{\mathrm{D}t}& = - \nabla p (\mathbf{x}, t) + \mu
  \nabla^2 \mathbf{u} (\mathbf{x}, t), & \mathbf{x} \in \Omega,  \\
  \nabla \cdot \mathbf{u} (\mathbf{x}, t) & = 0, & \mathbf{x} \in
  \Omega, \\
  \llbracket \mathbf{u} (\mathbf{\boldsymbol{\chi}} (\mathbf{X}, t), t) \rrbracket &
  = 0,  &\mathbf{X} \in \Gamma_0,  \\
  \llbracket p (\mathbf{\boldsymbol{\chi}} (\mathbf{X}, t), t) \rrbracket & = -\jmath^{-
  1} (\mathbf{X}, t)\, \mathbf{F} (\mathbf{X}, t) \cdot \mathbf{n} (\mathbf{X}, t), & \mathbf{X} \in \Gamma_0, \label{jump_condition_first_order_pressure}\\
  \mu\left\llbracket  \frac{\partial \mathbf{u} (\mathbf{\boldsymbol{\chi}}
  (\mathbf{X}, t), t)}{\partial x_i} \right\rrbracket &
  =   (\mathbb{I}- \mathbf{n} (\mathbf{X}, t)
  \mathbf{n} (\mathbf{X}, t)\tran) \jmath^{- 1} (\mathbf{X}, t)\, \mathbf{F} (\mathbf{X}, t) n_i  (\mathbf{X}, t), &
  \mathbf{X} \in \Gamma_0, \label{jump_condition_first_order_velocity}
  \\  \mathbf{U} (\mathbf{X}, t)  &
  = \mathbf{u} (\boldsymbol{\chi}(\mathbf{X}, t), t),  &\mathbf{X} \in \Gamma_0, 
  \\ \frac{\partial \boldsymbol{\chi(\mathbf{X},t)}}{\partial t}     &
  = \mathbf{U}(\mathbf{X},t),  &\mathbf{X} \in \Gamma_0, \label{eqn:interface_condition_7}
\end{align}
in which $\rho$, $\mu$, $\mathbf{u} (\mathbf{x}, t)$, $p
(\mathbf{x}, t)$, and  are the fluid's mass
density, dynamic viscosity, velocity, and pressure. 
$\mathbf{F} (\mathbf{X}, t)$ and $\mathbf{U} (\mathbf{X}, t)$ are the interfacial force and velocity, respectively, expressed using Lagrangian material coordinates, and $\mathbf{n} (\mathbf{X}, t)$ is the surface's unit normal vector in the interface's current configuration.
 $\jmath^{- 1} =
\frac{{\mathrm{d}a}}{{\mathrm{d}A}}$ relates surface area along the interface in the current $(\mathrm{d}a)$ and reference $(\mathrm{d}A)$ configurations. 
\iffalse
We also investigate fluid structure interactions in which the structure is elastically deformable.
In these particular studies, the structural dynamics are modeled by
\begin{align}
  \rho_\text{s} \frac{\partial^2 \boldsymbol{\xi}}{\partial t^2} (\mathbf{X}, t)  &
  = \nabla_\mathbf{X} \cdot \mathbb{P}^\text{s}(\mathbf{X}, t),  &\mathbf{X} \in \Gamma_0^\text{s}, 
  \\ J (\mathbf{X}, t) & = 1   &
  \mathbf{X} \in \Gamma_0^\text{s}, 
\end{align}
in which $J = \text{det}(\mathbb{F})$ is the volumetric Jacobian determinant for structural deformations, $\rho_\text{s}$ is the structural mass density in the reference configuration, and $\mathbb{F}(\mathbf{X}, t) = \frac{\partial \boldsymbol{\xi}}{\partial \mathbf{X}}$ is the deformation gradient tensor. $\mathbb{P}^\text{s}(\mathbf{X}, t)$ is the first Piola-Kirchoff stress tensor, which is defined for an incompressible hyperelastic material by 
\begin{align}
\mathbb{P}^\text{s}(\mathbf{X}, t) &= \frac{\partial \psi}{\partial \mathbb{F} (\mathbf{X}, t)} - \Phi(\mathbf{X}, t) \mathbb{F}^{-T}(\mathbf{X}, t)&\mathbf{X} \in \Gamma_0^\text{s}
\end{align}
in which $\psi(\mathbf{X}, t)$ is the strain-energy functional and $\Phi(\mathbf{X}, t)$ is the hydrostatic pressure that enforces the structure's incompressibility. 
\fi
The jump conditions express the discontinuity in fluid traction generated by the singular force along the interface. 
For a general function $\phi(\mathbf{x},t)$ and position $\mathbf{x}= \mathbf{\boldsymbol{\chi}}
  (\mathbf{X}, t)$ along the interface, the jump condition is
\begin{equation}
  \llbracket \phi (\mathbf{x}, t) \rrbracket \equiv \phi^+ (\mathbf{x}, t) - \phi^- (\mathbf{x}, t) \equiv \lim_{\epsilon
  \rightarrow 0^+} \phi (\mathbf{x} + \epsilon \mathbf{n}
  (\mathbf{X}, t), t) - \lim_{\varepsilon \rightarrow 0^-} \phi
  (\mathbf{x} - \varepsilon \mathbf{n} (\mathbf{X}, t), t).
\end{equation}
Equations (3)--(5) were derived in previous work by  Peskin and Printz\cite{PESKIN199333}, Lai and Li\cite{li_immersed_2001}, and Xu and Wang\cite{xu_systematic_2006}. Peskin and Printz showed that the delta function formulation used in the immersed boundary method, which directly applies singular forcing in the momentum equation (1), reduces to these jump conditions across the interface. Thus, our continuous equations of motion do not explicitly contain the Eulerian forcing. The effects of the jump conditions, however, arise in the discrete equations of motion through correction terms.

Our numerical tests focus on cases with stationary boundaries for which $\jmath^{- 1} = 1$. We view the boundary as having a \textit{prescribed} configuration $\boldsymbol{\xi} (\mathbf{X},t)$ and an \textit{actual} configuration $\mathbf{\boldsymbol{\chi}} (\mathbf{X}, t)$. 
\iffalse
In the case of elastic body motion, $\mathbf{\boldsymbol{\xi}} (\mathbf{X}, t))$ is determined by solving for the hyperelastic motion \cite{KOLAHDOUZ2021110442}.
\fi
Let $\mathbf{V}(\mathbf{X}, t) = \frac{\partial\boldsymbol{\xi} (\mathbf{X}, t)}{\partial t}$ be the velocity of the prescribed configuration. These two configurations are connected through penalty forces that are of the form

\begin{equation}
    \mathbf{F} (\mathbf{X}, t)  =  \kappa (\boldsymbol{\xi} (\mathbf{X},
    t) - \mathbf{\boldsymbol{\chi}} (\mathbf{X}, t)) + \eta \left( \mathbf{V}(\mathbf{X},
    t) - \mathbf{U}
    (\mathbf{X}, t) \right).
\end{equation}
Here, $\kappa$ and $\eta$ are the penalty spring stiffness and damping coefficients. As $\kappa \rightarrow \infty$, this constraint exactly
imposes the motion\cite{IMFSI-griffith}. In the simplest case, in which the interface is stationary, such that $\boldsymbol{\chi}(\mathbf{X},t) = \mathbf{X}$ and $\mathbf{V}(\mathbf{X},t) = \mathbf{0}$,
the simplified force model becomes $\mathbf{F} (\mathbf{X}, t) = \kappa (\mathbf{X} - \mathbf{\boldsymbol{\chi}} (\mathbf{X}, t)) - \eta  \mathbf{U}
    (\mathbf{X}, t) $.

\section{Discrete Equations of Motion}

This section describes the spatial and temporal discretizations of the equations of motion. To simplify
the notation, the numerical scheme is presented in two spatial dimensions. Like the original discrete IIM~\cite{kolahdouz_immersed_2020}, the extension of the method
to three spatial dimensions is straightforward.

\subsection{Interface Representation}
\iffalse
We consider a family of triangulation $\{\Tau_h \}$ for $\Gamma_0$. The domain formed by $\Tau_h $ is denoted by $\Gamma_{h,0}=\cup_{T\in \{\Tau_h \}}\overline{T}$. We use a standard $H^1(\Gamma_{h,0})$-conforming finite element space to represent interfacial Lagrangian variables in the system.  We denote the finite element space as $\mathscr{V}_h$. Consequently, in the discrete representation of the Equations~(\ref{jump_condition_first_order_pressure})--(\ref{eqn:interface_condition_7}), the Lagrangian variables $\bF,  \bU, \bchi$ are represented in the finite element space $[\mathscr{V}_h]^2$. For this study, we consider $\mathscr{V}_h$ to be standard $\bP^1$ finite element space. 
\fi
We use a finite element representation of the immersed interface with
triangulation $\mathcal{T}_h$ of $\Gamma_0$, the reference configuration.
Consider elements $U_i$ such that $\mathcal{T}_h = \bigcup_i U_i$, with $i$ indexing
the mesh elements. The nodes of the mesh elements are $\{
\mathbf{X}_j \}_{j = 1}^M$ and have corresponding nodal (Lagrangian) basis functions $\{ \psi_j (\mathbf{X})
\}_{j = 1}^M$. Herein, we consider piecewise linear interface representations, in which $\psi_j$ is a piecewise linear Lagrange polynomial, although extending the interface representation to other piecewice polynomials is straightforward. The nodal basis functions $\psi_j(\mathbf{X})$ are continuous across elements but are not differentiable across element boundaries. The current location of the interfacial
nodes at time $t$ are $\{ \mathbf{\boldsymbol{\chi}}_j (t) \}_{j = 1}^M$. In the finite
element space dictated by the subspace $S_h = \text{span} \{ \psi_j
(\mathbf{X}) \}_{j = 1}^M$, the configuration of the interface is
$ \mathbf{\boldsymbol{\chi}}_h (\mathbf{X}, t) = \sum_{j = 1}^M \mathbf{\boldsymbol{\chi}}_j (t)
  \psi_j (\mathbf{X})$.

\begin{figure}[t]
\center{}

\tikzset{every picture/.style={line width=0.75pt}} %set default line width to 0.75pt        

\begin{tikzpicture}[x=0.75pt,y=0.75pt,yscale=-1,xscale=1]
%uncomment if require: \path (0,300); %set diagram left start at 0, and has height of 300

%Shape: Grid [id:dp6945417026347169] 
\draw  [draw opacity=0] (146,47) -- (385.82,47) -- (385.82,263.82) -- (146,263.82) -- cycle ; \draw  [color={rgb, 255:red, 155; green, 155; blue, 155 }  ,draw opacity=1 ] (146,47) -- (146,263.82)(166,47) -- (166,263.82)(186,47) -- (186,263.82)(206,47) -- (206,263.82)(226,47) -- (226,263.82)(246,47) -- (246,263.82)(266,47) -- (266,263.82)(286,47) -- (286,263.82)(306,47) -- (306,263.82)(326,47) -- (326,263.82)(346,47) -- (346,263.82)(366,47) -- (366,263.82) ; \draw  [color={rgb, 255:red, 155; green, 155; blue, 155 }  ,draw opacity=1 ] (146,47) -- (385.82,47)(146,67) -- (385.82,67)(146,87) -- (385.82,87)(146,107) -- (385.82,107)(146,127) -- (385.82,127)(146,147) -- (385.82,147)(146,167) -- (385.82,167)(146,187) -- (385.82,187)(146,207) -- (385.82,207)(146,227) -- (385.82,227)(146,247) -- (385.82,247) ; \draw  [color={rgb, 255:red, 155; green, 155; blue, 155 }  ,draw opacity=1 ]  ;
%Shape: Rectangle [id:dp7780339158009218] 
\draw  [line width=1.5]  (146,47) -- (385.82,47) -- (385.82,263.82) -- (146,263.82) -- cycle ;
%Shape: Ellipse [id:dp6287542850135597] 
\draw  [line width=1.5]  (193.57,162.82) .. controls (201.97,128.58) and (234.65,100.82) .. (266.57,100.82) .. controls (298.49,100.82) and (317.56,128.58) .. (309.16,162.82) .. controls (300.76,197.06) and (268.07,224.82) .. (236.16,224.82) .. controls (204.24,224.82) and (185.17,197.06) .. (193.57,162.82) -- cycle ;
%Shape: Circle [id:dp6435520178922365] 
\draw  [fill={rgb, 255:red, 0; green, 0; blue, 0 }  ,fill opacity=1 ] (199,141) .. controls (199,138.79) and (200.79,137) .. (203,137) .. controls (205.21,137) and (207,138.79) .. (207,141) .. controls (207,143.21) and (205.21,145) .. (203,145) .. controls (200.79,145) and (199,143.21) .. (199,141) -- cycle ;
%Shape: Circle [id:dp1182976395827674] 
\draw  [fill={rgb, 255:red, 0; green, 0; blue, 0 }  ,fill opacity=1 ] (211,125) .. controls (211,122.79) and (212.79,121) .. (215,121) .. controls (217.21,121) and (219,122.79) .. (219,125) .. controls (219,127.21) and (217.21,129) .. (215,129) .. controls (212.79,129) and (211,127.21) .. (211,125) -- cycle ;
%Shape: Circle [id:dp24954158249256808] 
\draw  [fill={rgb, 255:red, 0; green, 0; blue, 0 }  ,fill opacity=1 ] (227,113) .. controls (227,110.79) and (228.79,109) .. (231,109) .. controls (233.21,109) and (235,110.79) .. (235,113) .. controls (235,115.21) and (233.21,117) .. (231,117) .. controls (228.79,117) and (227,115.21) .. (227,113) -- cycle ;
%Shape: Circle [id:dp6691301712755652] 
\draw  [fill={rgb, 255:red, 0; green, 0; blue, 0 }  ,fill opacity=1 ] (190.57,160.82) .. controls (190.57,158.61) and (192.36,156.82) .. (194.57,156.82) .. controls (196.78,156.82) and (198.57,158.61) .. (198.57,160.82) .. controls (198.57,163.03) and (196.78,164.82) .. (194.57,164.82) .. controls (192.36,164.82) and (190.57,163.03) .. (190.57,160.82) -- cycle ;
%Shape: Circle [id:dp08221558192552403] 
\draw  [fill={rgb, 255:red, 0; green, 0; blue, 0 }  ,fill opacity=1 ] (187.57,181.82) .. controls (187.57,179.61) and (189.36,177.82) .. (191.57,177.82) .. controls (193.78,177.82) and (195.57,179.61) .. (195.57,181.82) .. controls (195.57,184.03) and (193.78,185.82) .. (191.57,185.82) .. controls (189.36,185.82) and (187.57,184.03) .. (187.57,181.82) -- cycle ;
%Shape: Circle [id:dp5592914045234992] 
\draw  [fill={rgb, 255:red, 0; green, 0; blue, 0 }  ,fill opacity=1 ] (193.57,202.82) .. controls (193.57,200.61) and (195.36,198.82) .. (197.57,198.82) .. controls (199.78,198.82) and (201.57,200.61) .. (201.57,202.82) .. controls (201.57,205.03) and (199.78,206.82) .. (197.57,206.82) .. controls (195.36,206.82) and (193.57,205.03) .. (193.57,202.82) -- cycle ;
%Shape: Circle [id:dp17965171926094858] 
\draw  [fill={rgb, 255:red, 0; green, 0; blue, 0 }  ,fill opacity=1 ] (207.57,217.82) .. controls (207.57,215.61) and (209.36,213.82) .. (211.57,213.82) .. controls (213.78,213.82) and (215.57,215.61) .. (215.57,217.82) .. controls (215.57,220.03) and (213.78,221.82) .. (211.57,221.82) .. controls (209.36,221.82) and (207.57,220.03) .. (207.57,217.82) -- cycle ;
%Shape: Circle [id:dp6967272671844978] 
\draw  [fill={rgb, 255:red, 0; green, 0; blue, 0 }  ,fill opacity=1 ] (227.16,223.82) .. controls (227.16,221.61) and (228.95,219.82) .. (231.16,219.82) .. controls (233.37,219.82) and (235.16,221.61) .. (235.16,223.82) .. controls (235.16,226.03) and (233.37,227.82) .. (231.16,227.82) .. controls (228.95,227.82) and (227.16,226.03) .. (227.16,223.82) -- cycle ;
%Shape: Circle [id:dp48280937534422097] 
\draw  [fill={rgb, 255:red, 0; green, 0; blue, 0 }  ,fill opacity=1 ] (248.57,221.82) .. controls (248.57,219.61) and (250.36,217.82) .. (252.57,217.82) .. controls (254.78,217.82) and (256.57,219.61) .. (256.57,221.82) .. controls (256.57,224.03) and (254.78,225.82) .. (252.57,225.82) .. controls (250.36,225.82) and (248.57,224.03) .. (248.57,221.82) -- cycle ;
%Shape: Circle [id:dp19185535442552648] 
\draw  [fill={rgb, 255:red, 0; green, 0; blue, 0 }  ,fill opacity=1 ] (268.57,213.82) .. controls (268.57,211.61) and (270.36,209.82) .. (272.57,209.82) .. controls (274.78,209.82) and (276.57,211.61) .. (276.57,213.82) .. controls (276.57,216.03) and (274.78,217.82) .. (272.57,217.82) .. controls (270.36,217.82) and (268.57,216.03) .. (268.57,213.82) -- cycle ;
%Shape: Circle [id:dp35547102358076943] 
\draw  [fill={rgb, 255:red, 0; green, 0; blue, 0 }  ,fill opacity=1 ] (284.57,199.82) .. controls (284.57,197.61) and (286.36,195.82) .. (288.57,195.82) .. controls (290.78,195.82) and (292.57,197.61) .. (292.57,199.82) .. controls (292.57,202.03) and (290.78,203.82) .. (288.57,203.82) .. controls (286.36,203.82) and (284.57,202.03) .. (284.57,199.82) -- cycle ;
%Shape: Circle [id:dp003975365735318204] 
\draw  [fill={rgb, 255:red, 0; green, 0; blue, 0 }  ,fill opacity=1 ] (297.57,182.82) .. controls (297.57,180.61) and (299.36,178.82) .. (301.57,178.82) .. controls (303.78,178.82) and (305.57,180.61) .. (305.57,182.82) .. controls (305.57,185.03) and (303.78,186.82) .. (301.57,186.82) .. controls (299.36,186.82) and (297.57,185.03) .. (297.57,182.82) -- cycle ;
%Shape: Circle [id:dp7901041927816319] 
\draw  [fill={rgb, 255:red, 0; green, 0; blue, 0 }  ,fill opacity=1 ] (306,161) .. controls (306,158.79) and (307.79,157) .. (310,157) .. controls (312.21,157) and (314,158.79) .. (314,161) .. controls (314,163.21) and (312.21,165) .. (310,165) .. controls (307.79,165) and (306,163.21) .. (306,161) -- cycle ;
%Shape: Circle [id:dp13431417139092727] 
\draw  [fill={rgb, 255:red, 0; green, 0; blue, 0 }  ,fill opacity=1 ] (247,103) .. controls (247,100.79) and (248.79,99) .. (251,99) .. controls (253.21,99) and (255,100.79) .. (255,103) .. controls (255,105.21) and (253.21,107) .. (251,107) .. controls (248.79,107) and (247,105.21) .. (247,103) -- cycle ;
%Shape: Circle [id:dp5259625565792072] 
\draw  [fill={rgb, 255:red, 0; green, 0; blue, 0 }  ,fill opacity=1 ] (268.57,101.82) .. controls (268.57,99.61) and (270.36,97.82) .. (272.57,97.82) .. controls (274.78,97.82) and (276.57,99.61) .. (276.57,101.82) .. controls (276.57,104.03) and (274.78,105.82) .. (272.57,105.82) .. controls (270.36,105.82) and (268.57,104.03) .. (268.57,101.82) -- cycle ;
%Shape: Circle [id:dp8193151818330852] 
\draw  [fill={rgb, 255:red, 0; green, 0; blue, 0 }  ,fill opacity=1 ] (287,108) .. controls (287,105.79) and (288.79,104) .. (291,104) .. controls (293.21,104) and (295,105.79) .. (295,108) .. controls (295,110.21) and (293.21,112) .. (291,112) .. controls (288.79,112) and (287,110.21) .. (287,108) -- cycle ;
%Shape: Circle [id:dp41699986337944006] 
\draw  [fill={rgb, 255:red, 0; green, 0; blue, 0 }  ,fill opacity=1 ] (301,121) .. controls (301,118.79) and (302.79,117) .. (305,117) .. controls (307.21,117) and (309,118.79) .. (309,121) .. controls (309,123.21) and (307.21,125) .. (305,125) .. controls (302.79,125) and (301,123.21) .. (301,121) -- cycle ;
%Shape: Circle [id:dp5562173996633977] 
\draw  [fill={rgb, 255:red, 0; green, 0; blue, 0 }  ,fill opacity=1 ] (306.57,139.82) .. controls (306.57,137.61) and (308.36,135.82) .. (310.57,135.82) .. controls (312.78,135.82) and (314.57,137.61) .. (314.57,139.82) .. controls (314.57,142.03) and (312.78,143.82) .. (310.57,143.82) .. controls (308.36,143.82) and (306.57,142.03) .. (306.57,139.82) -- cycle ;
%Straight Lines [id:da9187663397972285] 
\draw    (337.82,205.82) -- (309.1,187.44) ;
\draw [shift={(306.57,185.82)}, rotate = 32.62] [fill={rgb, 255:red, 0; green, 0; blue, 0 }  ][line width=0.08]  [draw opacity=0] (8.93,-4.29) -- (0,0) -- (8.93,4.29) -- cycle    ;
%Curve Lines [id:da6256813221347833] 
\draw    (229,72) .. controls (219.38,74.66) and (225.1,89.99) .. (235.89,102.64) ;
\draw [shift={(237.82,104.82)}, rotate = 227.29] [fill={rgb, 255:red, 0; green, 0; blue, 0 }  ][line width=0.08]  [draw opacity=0] (8.93,-4.29) -- (0,0) -- (8.93,4.29) -- cycle    ;

% Text Node
\draw (155.36,58.4) node [anchor=north west][inner sep=0.75pt]  [font=\large]  {$\Omega _{t}^{+}$};
% Text Node
\draw (236.36,144.4) node [anchor=north west][inner sep=0.75pt]  [font=\large]  {$\Omega _{t}^{-}$};
% Text Node
\draw (312,207.4) node [anchor=north west][inner sep=0.75pt]  [font=\large]  {$\chi (\mathbf{X} ,t)$};
% Text Node
\draw (231,54.4) node [anchor=north west][inner sep=0.75pt]  [font=\normalsize]  {$ \begin{array}{l}
\Gamma _{t} =\ \overline{\Omega _{t}^{+}}\bigcap \overline{\Omega _{t}^{-}}\\
\end{array}$};

\end{tikzpicture}

  \caption{A Lagrangian coordinate system represents represents the interface, $\Gamma_t$. At time $t$, the location of the interface in Eulerian coordinates is $\boldsymbol{\chi}(\mathbf{X},t)$.}
\end{figure}
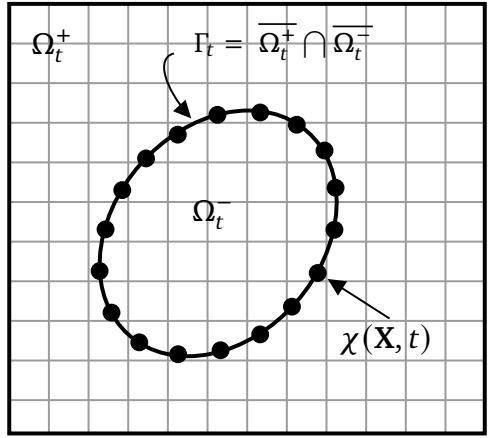

\subsection{Finite Difference Approximation}

We use a staggered-grid discretization of the incompressible Navier-Stokes equations. The discretization approximates 
pressure at cell centers, and velocity and forcing terms at the center of cell edges (in two spatial dimensions) or cell faces (in three spatial dimensions)\cite{griffith_volume_2012,griffith_accurate_2009}. Our computations use an isotropic grid, such that $h=\Delta x = \Delta y$. We use second-order accurate finite
difference stencils, and define the discrete divergence of the velocity
$\mathbf{D} \cdot \mathbf{u}$ on cell centers, and the discrete
pressure gradient $\mathbf{G} p$, and discrete Laplacian $L
\mathbf{u} $ on cell edges (or, in three spatial dimensions, faces). Some numerical experiments detailed herin use adaptive mesh refinement{\cite{griffith_immersed_2012}} to
increase computational efficiency. The mesh factor ratio,
\begin{align*}
    M_{\text{fac}}=\frac{\text{Lagrangian element size}}{\text{Eulerian grid size}},
\end{align*} 
describes the relative resolutions of the Lagrangian and Eulerian meshes.

\subsection{Force Spreading and Velocity Interpolation}

In modern IIMs, jump conditions are imposed through correction terms that appear on the right-hand side of the \textit{discretized} incompressible Navier--Stokes equations. Such numerical methods can be described and implemented in a way that is analogous to the conventional IB method, except that they use different operators to connect the discretized Lagrangian and Eulerian variables.
Formulations that target only the lowest order jump conditions, as considered herein, include corrections only in the discretized momentum equation.
Consequently, these corrections take the form of additional body forces that are concentrated along the immersed boundary.
Because the jump conditions are determined from the configuration of the boundary and the interfacial forces, it is convenient to introduce a force spreading operator, 
$\boldsymbol{\mathcal{S}}[\boldsymbol{\chi}]$, that relates the Lagrangian force
$\mathbf{F}_h (\mathbf{X}, t)  = \sum_{j = 1}^M  \mathbf{F}_j(t) \,\psi_j (\mathbf{X})$
to the Eulerian forces $\mathbf{f}$ on the Cartesian grid by summing all correction terms in the modified finite difference stencils involving the Lagrangian force. Correction terms, which are detailed in previous work \cite{kolahdouz_immersed_2020}, depend on the projected jump conditions, as computed in Section 3.4.

To derive these corrections, consider the case in which the interface intersects a grid line between two grid points at  $x=\alpha$, $x_i\leq \alpha\leq x_{i+1}$. For a twice differentiable piecewise function $v(x)$, we approximate its first and second derivatives as using the lowest order jump conditions:

\iffalse
\begin{align}
    v'(x_i)&=\frac{v_{i+1}-v_{i-1}}{2h}-\frac{1}{2h}\sum_{m=0}^2\frac{(h^+)^m}{m!}\jumpinverse{v^{(m)}}_{\alpha}+\mathcal{O}(h^2) \label{eqn:taylor_1}\\
      v''(x_i)&=\frac{v_{i+1}-2v_i+v_{i-1}}{h^2}-\frac{1}{h^2}\sum_{m=0}^2\frac{(h^+)^m}{m!}\jumpinverse{v^{(m)}}_{\alpha}+\mathcal{O}(h^2)\label{eqn:taylor_4}
\end{align}
\fi

\begin{align}
    v'(x_{i+\frac{1}{2}})&=\frac{v_{i+1}-v_{i}}{h}-\frac{\jumpinverse{v}_{\alpha}}{h}, \label{eqn:taylor_1}\\
      v''(x_i)&=\frac{v_{i+1}-2v_i+v_{i-1}}{h^2}-\frac{h^+\jumpinverse{\frac{\partial v}{\partial x}}_{\alpha}}{h^2},\label{eqn:taylor_4}
\end{align}
in which $v_{i}=v(x_i)$, $h^+=x_{i+1}-\alpha$, $h^-=x_{i}-\alpha$ and $h$ is the Cartesian grid size. To construct the force spreading operator $\mathcal{S}[\boldsymbol{\chi} ](\bF)$, we substitute the interface jump condition Equations~(\ref{jump_condition_first_order_pressure})--(\ref{jump_condition_first_order_velocity}) into generalized Taylor expansion Equations~(\ref{eqn:taylor_1})--(\ref{eqn:taylor_4}), which are used to derive finite difference approximations of the differential operators in the momentum equation on stencils that intersect the interface.  The force terms are then moved to the right-hand side of the momentum equation, and we denote all such contributions by $\mathcal{S}[\boldsymbol{\chi}](\bF)$. For more details, we refer the reader to prior studies \cite{LeVeque97,LI2001822,xusheng06}.

The configuration of the Lagrangian mesh is updated by enforcing the no-slip no-penetration condition $\mathbf{U}_h = \mathcal{I}[{\boldsymbol{\chi},\mathbf{F}_\tau} ](\mathbf{u})$, in which $\mathcal{I}$ is the velocity interpolation operator that includes an additional correction term from Equation~(\ref{jump_condition_first_order_velocity}) to improve accuracy. To compute interpolated velocities, the interpolation operator requires only the jump conditions for the velocity gradient. Consequently, only the tangential component of the Lagrangian force, $\mathbf{F}_\tau$, is needed to define this operator.
For additional details, we refer the reader to  Tan et al.~\cite{tan2009}.

\subsection{Projected Jump Conditions}
\begin{figure}[t]
\centering
\hskip -.3cm
\includegraphics[width=.49\textwidth]{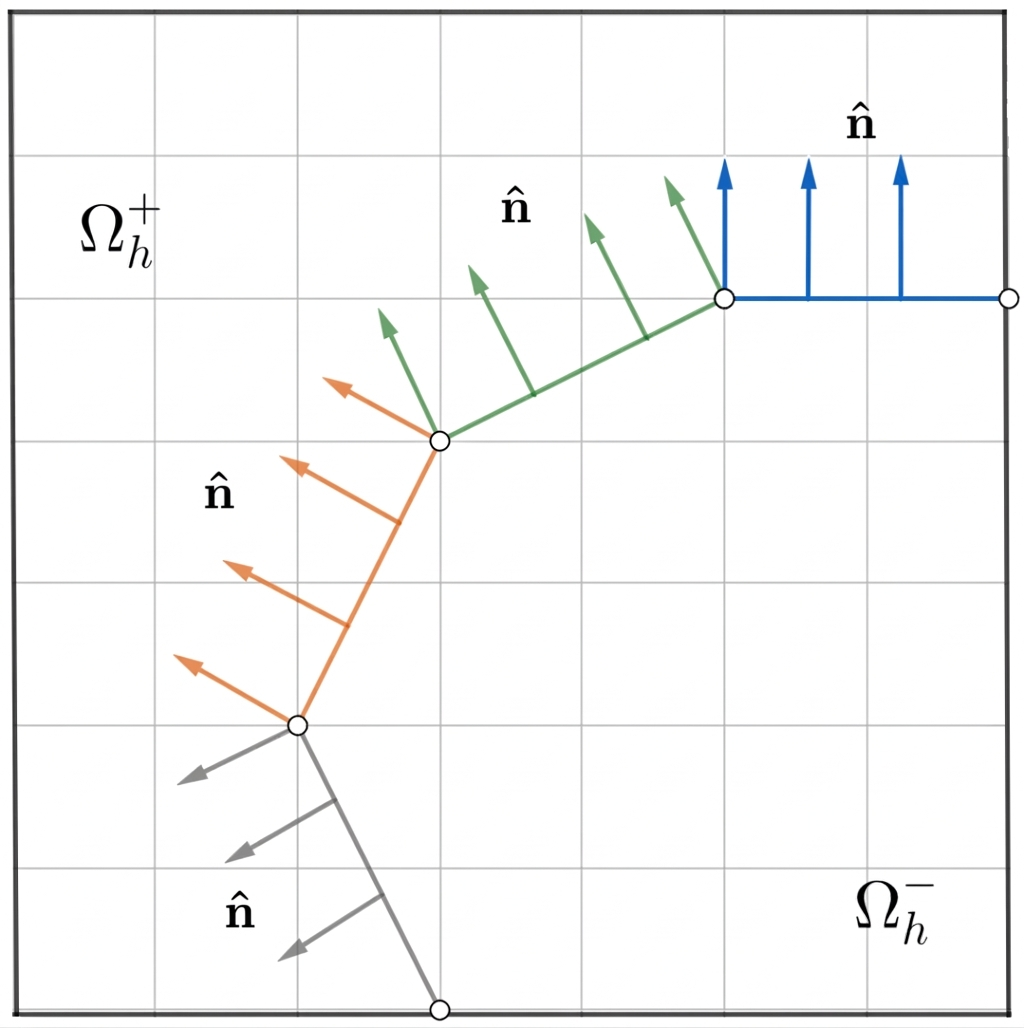}\quad
\includegraphics[width=.49\textwidth]{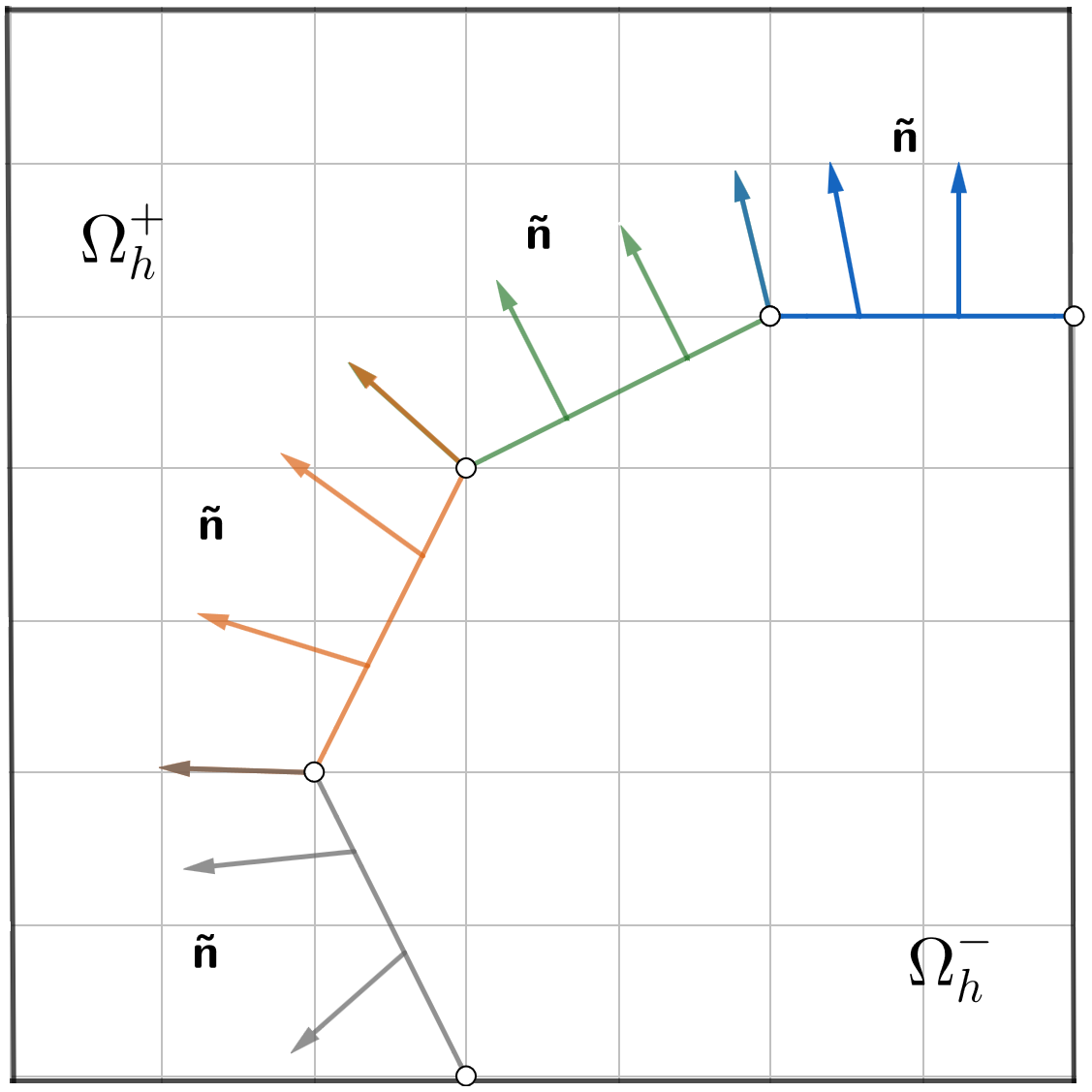}
  \caption{Two-dimensional schematic of a discrete immersed boundary. Normal vectors $\hat\bn$ are discontinuous at the junctions between neighboring elements (left). Globally continuous normal vector  field $\tilde{\bn}$  (right).}
  \label{fig:compare_normal}
\end{figure}

In our methodology, jump conditions must be evaluated at arbitrary positions along the interface to evaluate the forcing terms associated with stencils that cut the interface. 
Let $\hat \bn $ be the element-aligned outward normal on $\Gamma_{h,t} = \bchi(\Gamma_{h,0}, t)$, pointing from $\Omega_{h,t}^-$ to $\Omega_{h,t}^+$, denoted as $\hat \bn = (n^x, n^y)$. 
We are using both a local projection operator $\mathbf{P}(\mathbf{X})$ that operates on individual vectors, and a global projection operator $\boldsymbol{\mathcal{P}}_h$ that operates on fields.
Define the orthogonal projection onto the tangential plane as $\mathbf{P}(\mathbf{X}) = \mathbb{I} - \hat{\mathbf{n}}(\boldsymbol{\chi}(\mathbf{X}, t)) \hat{\mathbf{n}}\tran(\boldsymbol{\chi}(\mathbf{X}, t))$.
Because the normal vector field $\hat{\mathbf{n}}$ on $\Gamma_{h,t}$ is generally not continuous (see Figure~\ref{fig:compare_normal}), this leads to $\mathbf{F}_\tau$ and $F_\mathbf{n}$ being discontinuous in Equations~(\ref{jump_condition_first_order_pressure})--(\ref{jump_condition_first_order_velocity}). 
Because the interpolation operator $ \mathcal{I}[\boldsymbol{\chi},\mathbf{F}_\tau]$ contains a term based on Equations~(\ref{jump_condition_first_order_velocity}), if $\mathbf{F}_\tau$ is discontinuous, then the interpolated surface velocity $ \mathcal{I}[\boldsymbol{\chi},\mathbf{F}_\tau] (\mathbf{u})$ is also discontinuous.
To obtain a continuous representation of the interface velocity, previous work by Kolahdouz et al.\cite{kolahdouz_immersed_2020} proposed using a standard $L^2$ projection, $\boldsymbol{\mathcal{P}}_h$, to project both discontinuous pressure and shear stress jump conditions onto the continuous finite element space $\mathscr{V}_h$. The $L^2$ projection operator $\boldsymbol{\mathcal{P}}_h$ is defined in the classical sense:\\
For a given function $\bphi\in [L^2(\Gamma_{h,0})]^d$, we seek to find $\boldsymbol{\mathcal{P}}_h(\bphi)\in [\mathscr{V}_h]^d$ that satisfies,
\begin{align}
\int_{\Gamma_{h,0}}\boldsymbol{\mathcal{P}}_h(\bphi)\cdot \bpsi\ \mathrm{d}A=\int_{\Gamma_{h,0}}\bphi\cdot\bpsi\ \mathrm{d}A &\qquad \text{for any }  \bpsi\in [\mathscr{V}_h]^d.\label{eqn:mass_matrix}
\end{align}
 Our earlier work ~\cite{kolahdouz_immersed_2020} uses the following $L^2$-projected jump conditions in the discrete system:
\begin{align}
\left\llbracket\mu\frac{\partial \bu}{\partial x}\right\rrbracket =-\boldsymbol{\mathcal{P}}_h(\jmath^{-1}\bF_\tau n^x),\qquad \left\llbracket\mu\frac{\partial \bu}{\partial y}\right\rrbracket&=-\boldsymbol{\mathcal{P}}_h(\jmath^{-1}\bF_\tau n^y),\qquad \jumpinverse{p}= \boldsymbol{\mathcal{P}}_h(\jmath^{-1}F_\mathbf{n} ).\label{jump_condition_first_order_weak}
\end{align}
Equations~(\ref{jump_condition_first_order_weak}) are used in the construction of both the force spreading operator  $\mathcal{S}[\boldsymbol{\chi}](\bF)$ and the velocity interpolation operator $ \mathcal{I}[\boldsymbol{\chi},\mathbf{F}_t] (\mathbf{u})$, ensuring continuity of the interface velocity $\bU$. 

In previous work, each element's normal vector is used to compute the projected jump conditions. This study investigates the impact of alternative approaches that have been demonstrated to provide improved accuracy in computing surface curvature. Inspired by the vertex normal computation by  Chen et al.~\cite{CHEN2004447}, we use two methods to construct globally continuous normal vectors by computing vertex normals, $\mathbf{n}_i$, which are linearly interpolated across element interiors.

\subsection{Globally Continuous Normal Vectors}
This section introduces two straightforward methods for constructing the continuous normal and tangential surface vector fields. To do so, we first define the vertex normal, $\mathbf{n}_i$, because the discrete normal vector field, $\hat{\bn}(\mathbf{X}, t)$, which is determined from the $C^0$ geometry is discontinuous at vertices. For brevity, we will refer to this normal vector as a elemental-aligned normal. 

A straightforward way to construct $\mathbf{n}_i$ is to use the $L^2$ projection operator $\boldsymbol{\mathcal{P}}_h$ to project $\hat{\bn}=\hat{\bn}(\mathbf{X}, t)$ onto the $\mathbf{P}^1$ finite element basis.
We compute the $L^2$-projected vertex normal vectors $\mathbf{n}_i^\text{proj}$ via:
\begin{align}
    \mathbf{n}_i^\text{proj}=\frac{\boldsymbol{\mathcal{P}}_h(\hat{\mathbf{n}})}{\norm{\boldsymbol{\mathcal{P}}_h(\hat{\mathbf{n}})}_2}(\mathbf{X}_i).
\end{align}
Another way to compute the vertex normal vectors $\mathbf{n}_i$ is to use an inverse centroid-weighted (ICW) average of the neighboring elements that share vertex $i$. The weighted average of the neighboring elemental-aligned normals is:
\begin{equation}
 \mathbf{n}_i^\text{ICW}= \frac{ \sum_j \hat{\mathbf{n}}_j(\mathbf{X}^\text{C}_j, t) w_{i,j}
   }{\left\| \sum_j \hat{\mathbf{n}}_j (\mathbf{X}_j^\text{C}, t)w_{i,j}\right\|_2}, 
\end{equation}
in which $\hat{\mathbf{n}}_j$, $\mathbf{X}^\text{C}_j$, and $w_{i,j}$ are  the $j^\text{th}$ neighboring element's elemental-aligned normal, centroid location, and weight for the $i^\text{th}$ vertex. We compute  $w_{i,j} = \frac{1}{\ell_{i,j}}$, in which  $\ell_{i,j} = |\mathbf{X}_i-\mathbf{X}_j^\text{C}|$. This approach is similar to the vertex normal computations by Chen et al.~\cite{CHEN2004447}, in which the element's weight is computed using the inverse squared distance from an element's vertex to its centroid. 

To evaluate the unit normal vector at an element's interior for either definition of the nodal normal vectors, we use the $\mathbf{P}^1$ element's basis functions:
\begin{equation}
    \tilde{\mathbf{n}} (\mathbf{X}, t) = \frac{\sum_{j = 1}^M \mathbf{n}_j
  \psi_j(\mathbf{X})}{||\sum_{j = 1}^M \mathbf{n}_j
  \psi_j(\mathbf{X})||_2}.
\end{equation}

\subsection{Time Integration}

Each step begins with known values of $\mathbf{\boldsymbol{\chi}}^n$ and $\mathbf{u}^n$ at time
$t^n$, and $p^{n - \frac{1}{2}}$ at time $t^{n - \frac{1}{2}}$. The goal is to
compute $\mathbf{\boldsymbol{\chi}}^{n + 1}, \mathbf{u}^{n + 1}$, and $p^{n +
\frac{1}{2}}$. First, an initial prediction of the structure location at time $t^{n + 1}$ is determined by
\begin{equation}
  \frac{\widehat{\mathbf{\boldsymbol{\chi}}}^{n + 1} - \mathbf{\boldsymbol{\chi}}^n}{\Delta t} =
  \mathbf{U}^n (\mathbf{\boldsymbol{\chi}}^n) = \mathcal{I}[\boldsymbol{\chi}^n,\mathbf{F}^n] (\mathbf{u}^n)  .
\end{equation}
We can then approximate the structure location at time $t^{n + \frac{1}{2}}$ by
\begin{equation}
  \mathbf{\boldsymbol{\chi}}^{n + \frac{1}{2}} = \frac{\widehat{\mathbf{\boldsymbol{\chi}}}^{n + 1}
  + \mathbf{\boldsymbol{\chi}}^n}{2}.
\end{equation}
Next, we solve for $\mathbf{u}^{n + 1}$ and $p^{n +
\frac{1}{2}}$ in which $\mathbf{f}^{n + \frac{1}{2}} = \boldsymbol{\mathcal{S}}[\boldsymbol{\chi}^{n + \frac{1}{2}} ] \mathbf{F}(\boldsymbol{\chi}^{n + \frac{1}{2}},t^{n + \frac{1}{2}})$:
\begin{eqnarray}
  \rho \left( \frac{\mathbf{u}^{n + 1} - \mathbf{u}^n}{\Delta t} +
  \mathbf{A}^{n+\frac12} \right) & = & - \mathbf{G} p^{n +
  \frac{1}{2}} + \mu L \left( \frac{\mathbf{u}^{n + 1} +
  \mathbf{u}^n}{2} \right) + \mathbf{f}^{n + \frac{1}{2}}, \\
  \mathbf{D} \cdot \mathbf{u}^{n + 1} & = & 0,
\end{eqnarray}
in which the non-linear advection term, $\mathbf{A}^{n + \frac{1}{2}} = \frac{3}{2}\mathbf{A}^{n} - \frac{1}{2}\mathbf{A}^{n-1}$, is handled with the xsPPM7 variant\cite{xsPPM7} of the piecewise
parabolic method.{\cite{colella_piecewise-parabolic_1982}} 
$\mathbf{G}, \mathbf{D} \cdot\mbox{}$, and $L$ are the discrete
gradient, divergence, and Laplacian operators. The system of equations is
iteratively solved via the FGMRES
algorithm along with a preconditioner based on the projection method\cite{griffith_accurate_2009}. Last, we update the structure's location, $\mathbf{\boldsymbol{\chi}}^{n + 1}$, with
\begin{equation}
  \frac{\mathbf{\boldsymbol{\chi}}^{n + 1} - \mathbf{\boldsymbol{\chi}}^n}{\Delta t}  = 
  \mathbf{U}^{n + \frac{1}{2}} =  \boldsymbol{\mathcal{I}}[\boldsymbol{\chi}^{n+\frac{1}{2}},\mathbf{F}^{n+\frac{1}{2}}] \left(\frac{\mathbf{u}^n+ \mathbf{u}^{n+1}}{2}\right)  .
\end{equation}

\subsection{Software Implementation}

All computations for were completed through IBAMR \cite{ibamr}, which
utilizes parallel computing libraries and adaptive mesh refinement (AMR).
IBAMR uses other libraries for setting up meshes, fast linear algebra solvers,
and postprocessing, including SAMRAI \cite{samrai}, PETSc \cite{petsc,petsc2,petsc3}, \textit{hypre} \cite{hypre,hypre2}, and libMesh \cite{libmesh,libmesh2}.

\section{Numerical Experiments}

In this section, we investigate the effects of the different normal vector representations on inherently smooth geometries in both two and three spatial dimensions. 

We first examine flow past a circular cylinder as a benchmark case to examine the accuracy of projected normal vectors for external, bluff body flows. Next, we test our methods ability to enforce the no-penetration condition for our immersed interface with a pressurized cylinder in two dimensions.
We extend our tests to three dimensions, examining the accuracy for a pressure loaded cylinder, as well as Poiseuille flow.
We impose stationary boundaries using a penalty method described in Section 2. Because the penalty parameter is finite, generally there are discrepancies between the prescribed and actual positions of the interface. In numerical tests, we choose $\kappa$ to ensure that the discrepancy in interface configurations $\epsilon_{\mathbf{X}} = \| \mathbf{\xi} (\mathbf{X}, t) -
\mathbf{\boldsymbol{\chi}} (\mathbf{X}, t) \|$ satisfies $\epsilon_{\mathbf{X}} < \frac{h}{4}$. 
\

\subsection{Flow Past a Circular Cylinder}
This section considers flow past a stationary cylinder of diameter $D=1$ centered at the origin of the computational domain $\Omega = [- 15, 45] \times [-30, 30]$ with $L=60$. 
For $\mathbf{u} = (u, v)$, we use incoming flow velocity $\mathbf{} u = 1$, and $v = \cos \left( \pi \frac{y}{L} \right) e^{- 2 t}$
to initiate vortex shedding at a consistent time. We set $\rho = 1$ and $\mu = \frac{1}{\text{Re}}$, in which $\text{Re}$ is the Reynolds number. We use
$\Rey = \frac{\rho U D}{\mu}=200$ to ensure that we observe vortex
shedding.  On the coarsest level
of the hierarchical grid, $n_{\text{c}} = 32$ grid cells are used in both directions.
The effective number of grid cells on the $\ell^{\text{th}}$ level is
${n_{\ell}}  = 2^{\ell - 1} n_{\text{c}}$. For convergence studies, we fix $n_{\text{c}}$ and vary $\ell_{\text{max}}$, the maximum number of grid levels, from six to eight. For comparisons to prior works, we use eight levels. The time step size $\Delta t$ scales with $h_{\text{finest}}$ such that $\Delta t = \frac{
h_{\text{finest}}}{20}$, in which $h_{\text{finest}} = \frac{L}{{n_{\ell {\text{max}}}}}$. For this time step size and grid spacing, the Courant-Friedrichs-Lewy (CFL) number is approximately 0.04 once the model reaches periodic steady state. 
$\kappa$ and $\eta$ are scaled with the grid resolution on the finest level, such that $\kappa = \frac{C_\kappa}{h_{\text{finest}}} $ and  $\eta = \frac{C_\eta}{h_{\text{finest}}}$ with $C_\kappa = 3.413$ and $ C_\eta = 0.025 $ using a bisection method. The outflow
boundary uses zero normal and tangential traction, and the top and bottom boundaries use zero tangential traction and $v = 0$.
\begin{figure}[H]
    \centering

    % ==========================
    % ROW 1: Two Images Side-by-Side
    % ==========================
    
    % --- LEFT IMAGE (Flat Normals) ---
    \begin{minipage}[b]{0.49\textwidth}
        \centering
        % The Image
        \includegraphics[height=4cm]{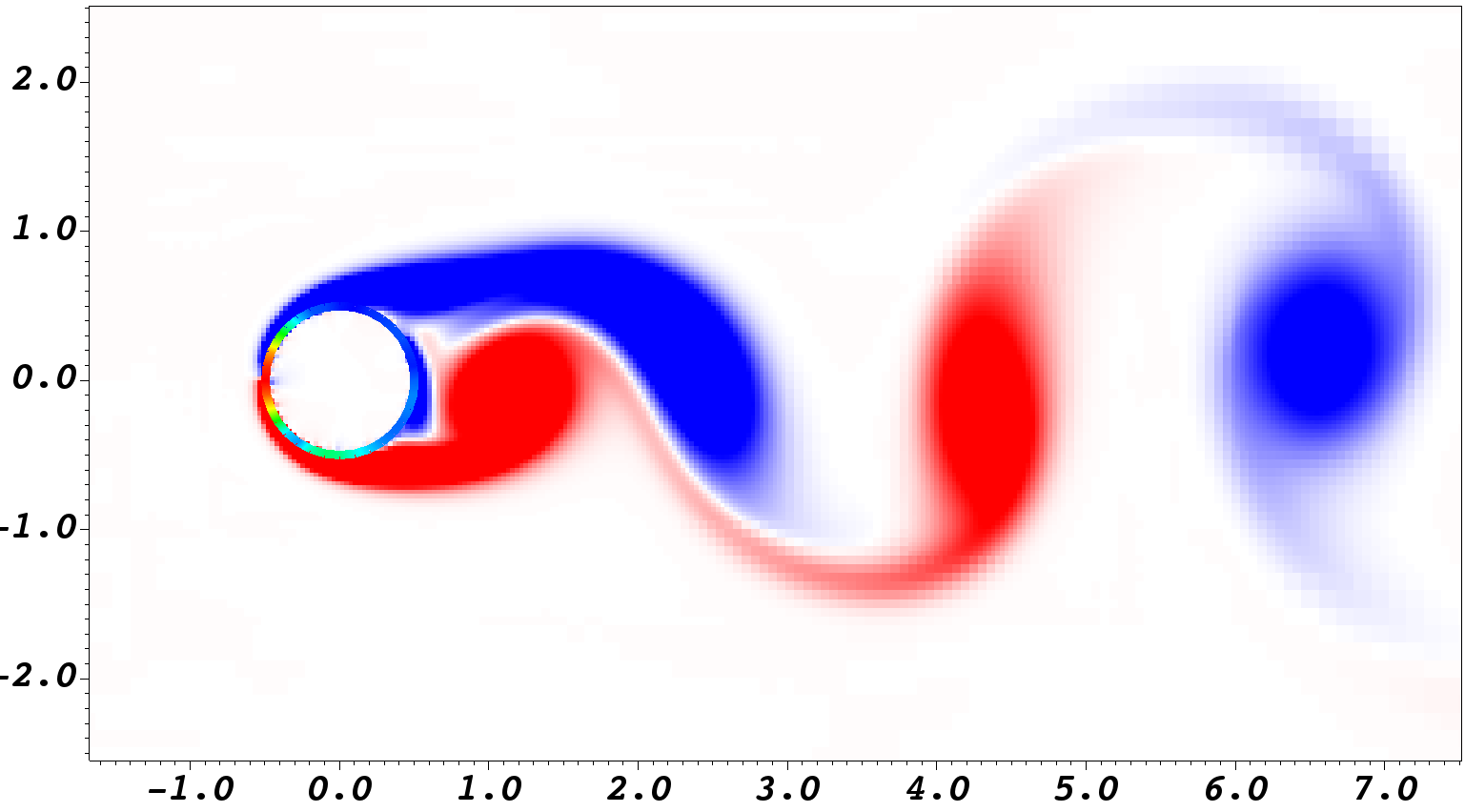}%
        % The Colorbar (Scaled to match image height)
        \resizebox{!}{4cm}{
            \tikzset{every picture/.style={line width=0.75pt}} 
            \begin{tikzpicture}[x=0.75pt,y=0.75pt,yscale=-1,xscale=1]
                \draw (56.38,216.48) node {\includegraphics[width=17.07pt,height=238.2pt]{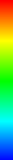}};
                \draw [color=black] (67.76,58.68) -- (75.76,58.68);
                \draw (67.76,375.28) -- (76.76,375.28);
                \draw (67.94,215.76) -- (76.76,215.76);
                \draw [color=black] (67.94,137.98) -- (76.76,137.98);
                \draw (67.94,292.54) -- (76.76,292.54);
                \draw (48,24) node [anchor=north west][font=\Large] {$\epsilon_{\mathbf{X}}$};
                \draw (72.82,47.42) node [anchor=north west][align=left] {9.9e-3};
                \draw (73.82,126.72) node [anchor=north west][align=left] {9.2e-3};
                \draw (75.82,204.02) node [anchor=north west][align=left] {8.5e-3};
                \draw (75.82,280.8) node [anchor=north west][align=left] {7.7e-3};
                \draw (77.82,360.58) node [anchor=north west][align=left] {7.0e-3};
            \end{tikzpicture}
        }
    \end{minipage}%
    \hfill % Pushes the images to the edges
    % --- RIGHT IMAGE (Nodal Normals) ---
    \begin{minipage}[b]{0.49\textwidth}
        \centering
        % The Image
        \includegraphics[height=4cm]{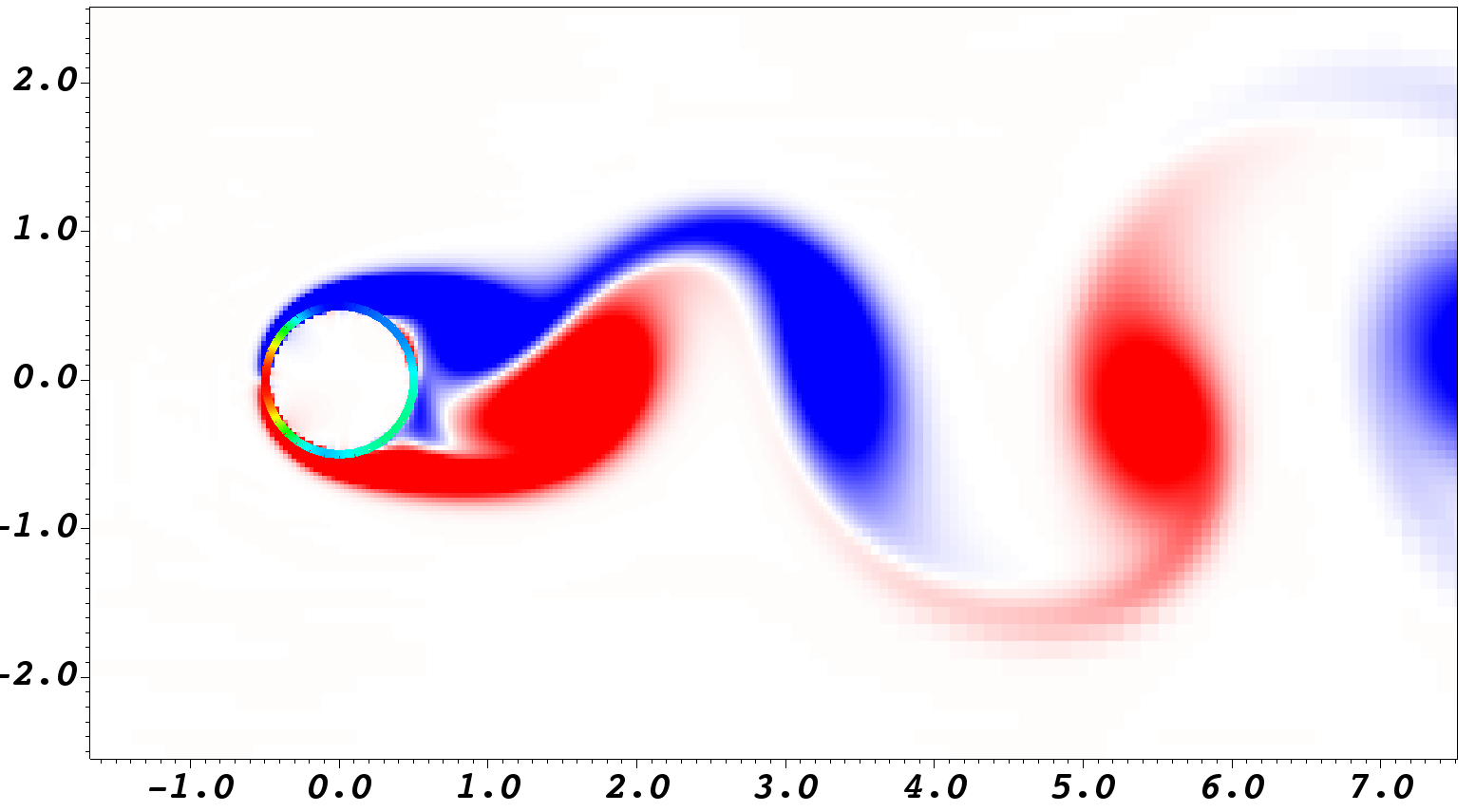}%
        % The Colorbar
        \resizebox{!}{4cm}{
            \tikzset{every picture/.style={line width=0.75pt}} %set default line width to 0.75pt        

\begin{tikzpicture}[x=0.75pt,y=0.75pt,yscale=-1,xscale=1]
%uncomment if require: \path (0,384); %set diagram left start at 0, and has height of 384

%Image [id:dp3783569225831278] 
\draw (56.38,216.48) node  {\includegraphics[width=17.07pt,height=238.2pt]{colorbar_displacement.png}};
%Straight Lines [id:da6763773634376069] 
\draw [color={rgb, 255:red, 0; green, 0; blue, 0 }  ,draw opacity=1 ]   (67.76,58.68) -- (75.76,58.68) ;
%Straight Lines [id:da1207005768568743] 
\draw    (67.76,375.28) -- (76.76,375.28) ;
%Straight Lines [id:da48914525610221027] 
\draw    (67.94,215.76) -- (76.76,215.76) ;
%Straight Lines [id:da7189166377438487] 
\draw [color={rgb, 255:red, 0; green, 0; blue, 0 }  ,draw opacity=1 ]   (67.94,137.98) -- (76.76,137.98) ;
%Straight Lines [id:da7702940702509451] 
\draw    (67.94,292.54) -- (76.76,292.54) ;

% Text Node
\draw (48,24) node [anchor=north west][inner sep=0.75pt]  [font=\Large]  {$\epsilon_{\mathbf{X}}$};
% Text Node
\draw (72.82,47.42) node [anchor=north west][inner sep=0.75pt]   [align=left] {2.1e-3};
% Text Node
\draw (73.82,126.72) node [anchor=north west][inner sep=0.75pt]   [align=left] {1.6e-3};
% Text Node
\draw (75.82,204.02) node [anchor=north west][inner sep=0.75pt]   [align=left] {1.1e-3};
% Text Node
\draw (75.82,280.8) node [anchor=north west][inner sep=0.75pt]   [align=left] {5.8e-4};
% Text Node
\draw (77.82,360.58) node [anchor=north west][inner sep=0.75pt]   [align=left] {5.8e-5};

\end{tikzpicture}
        }
    \end{minipage}

    \vspace{0.5cm} % Space between rows

    % ==========================
    % ROW 2: One Image Centered
    % ==========================
    
    % --- BOTTOM IMAGE (Consistent Normals) ---
    % Note: Made slightly larger (height=5cm) since it has more room
    \begin{minipage}[b]{0.8\textwidth} 
        \centering
        \resizebox{!}{5cm}{
            \tikzset{every picture/.style={line width=0.75pt}} 
            \begin{tikzpicture}[x=0.75pt,y=0.75pt,yscale=-1,xscale=1]
                \draw (56.38,216.48) node {\includegraphics[width=17.07pt,height=238.2pt]{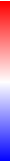}};
                \draw [color=black] (67.76,58.68) -- (75.76,58.68);
                \draw (67.76,375.28) -- (76.76,375.28);
                \draw (67.94,215.76) -- (76.76,215.76);
                \draw [color=black] (67.94,137.98) -- (76.76,137.98);
                \draw (67.94,292.54) -- (76.76,292.54);
                \draw (48,24) node [anchor=north west][font=\Large] {$\omega$};
                \draw (72.82,47.42) node [anchor=north west][align=left] {2.0};
                \draw (73.82,126.72) node [anchor=north west][align=left] {1.0};
                \draw (75.82,204.02) node [anchor=north west][align=left] {0.0};
                \draw (75.82,280.8) node [anchor=north west][align=left] {-2.0};
                \draw (77.82,360.58) node [anchor=north west][align=left] {-4.0};
            \end{tikzpicture}
        }
        % The Image
        \includegraphics[height=5cm]{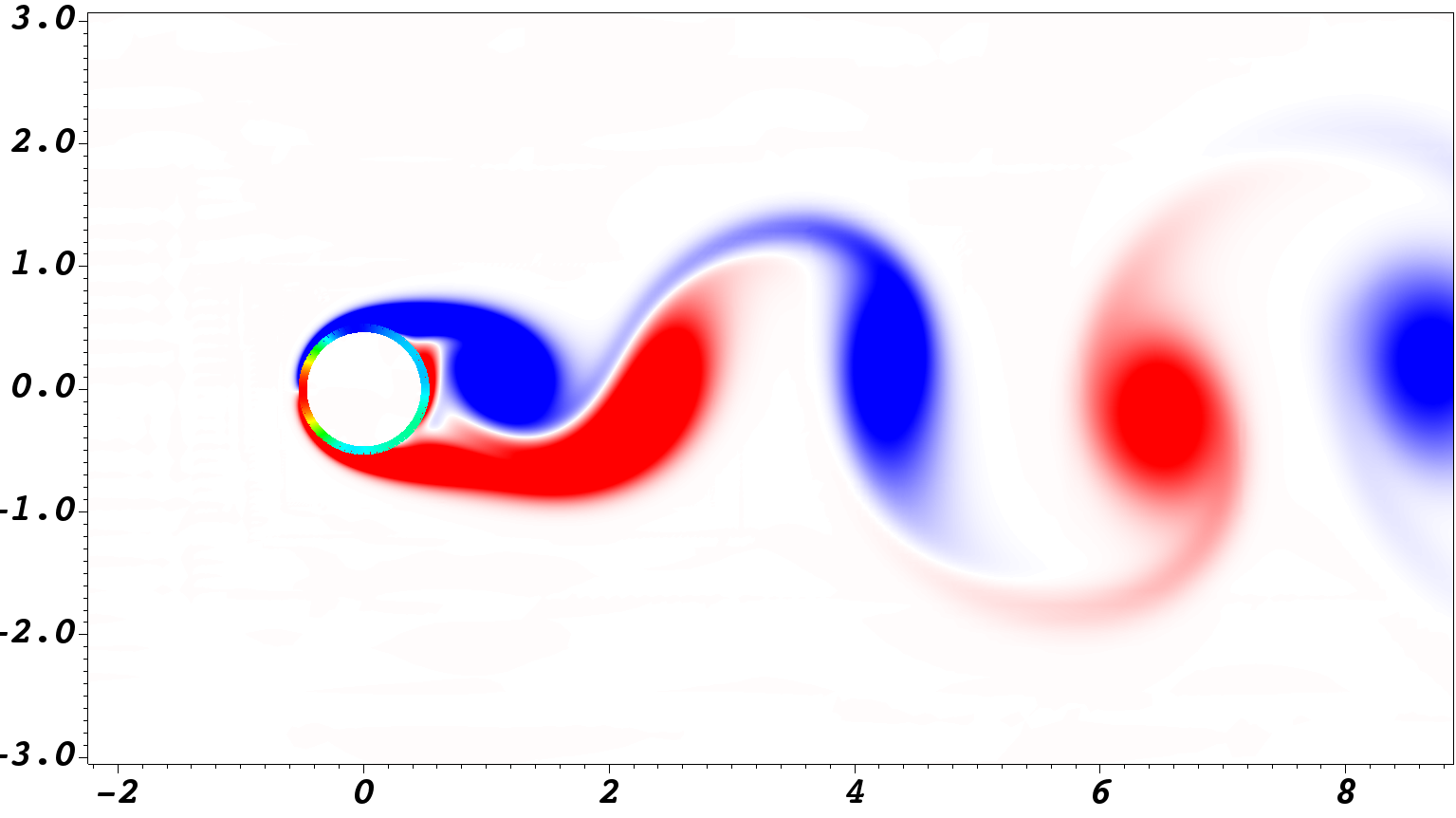}%
        % The Colorbar
        \resizebox{!}{5cm}{
            \tikzset{every picture/.style={line width=0.75pt}} 
            \begin{tikzpicture}[x=0.75pt,y=0.75pt,yscale=-1,xscale=1]
                \draw (56.38,216.48) node {\includegraphics[width=17.07pt,height=238.2pt]{colorbar_displacement.png}};
                \draw [color=black] (67.76,58.68) -- (75.76,58.68);
                \draw (67.76,375.28) -- (76.76,375.28);
                \draw (67.94,215.76) -- (76.76,215.76);
                \draw [color=black] (67.94,137.98) -- (76.76,137.98);
                \draw (67.94,292.54) -- (76.76,292.54);
                \draw (48,24) node [anchor=north west][font=\Large] {$\epsilon_{\mathbf{X}}$};
                \draw (72.82,47.42) node [anchor=north west][align=left] {7.2e-3};
                \draw (73.82,126.72) node [anchor=north west][align=left] {7.0e-3};
                \draw (75.82,204.02) node [anchor=north west][align=left] {6.8e-3};
                \draw (75.82,280.8) node [anchor=north west][align=left] {6.6e-3};
                \draw (77.82,360.58) node [anchor=north west][align=left] {6.5e-3};
            \end{tikzpicture}
        }
    \end{minipage}

    \caption{The top left, top right, and bottom panels show the vorticity field and $\varepsilon_{\mathbf{X}}$ computed using flat, inverse centroid-weighted, and projected normals with eight levels of grid refinement.}
\end{figure}

Figure 3 shows a snapshot of the vorticity field, $\omega$, at $t = 95$ for the three normal vector representations. Both flow fields show the expected vortex shedding and a max$(\varepsilon_{\mathbf{X}})$ on the same order of magnitude.

\begin{figure}[H]

  \centering
  \resizebox{200pt}{!}{\includegraphics{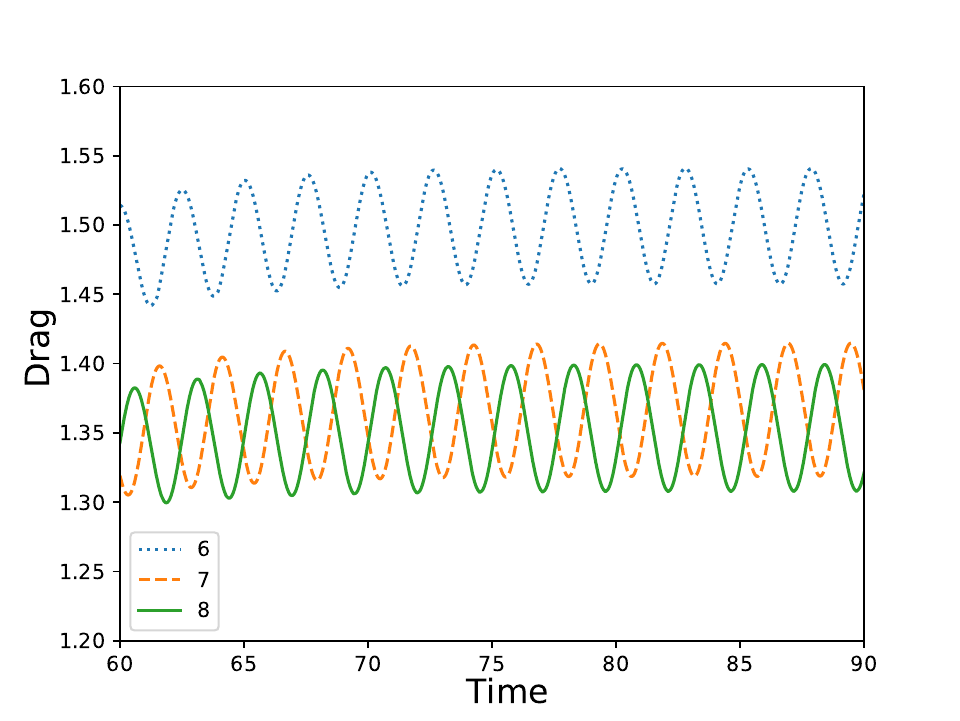}}
      \resizebox{200pt}{!}{\includegraphics{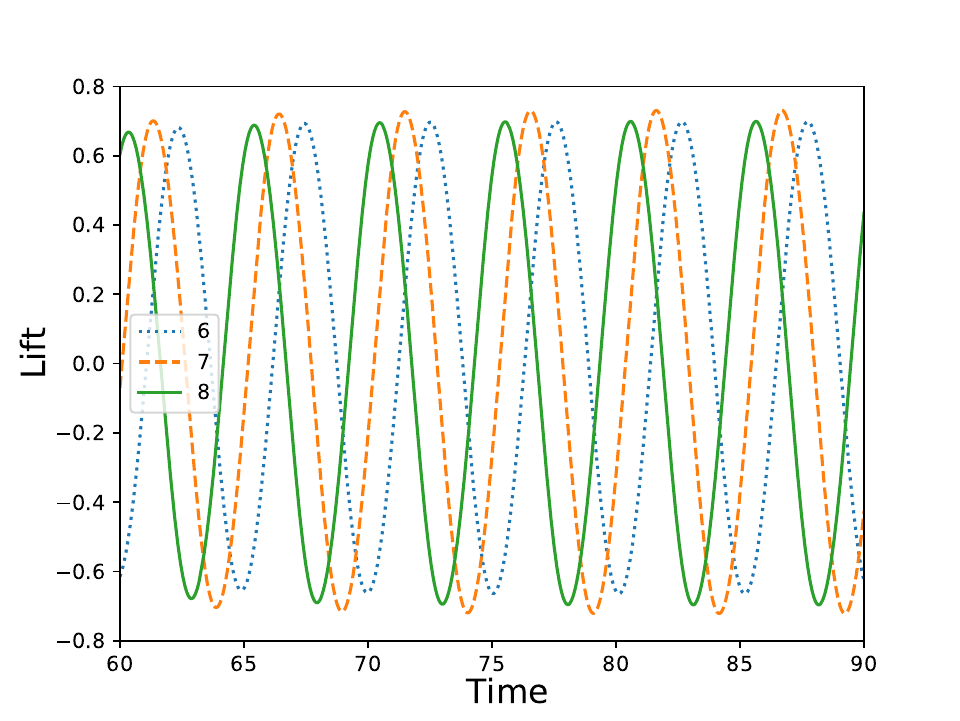}}
      
  \resizebox{200pt}{!}{\includegraphics{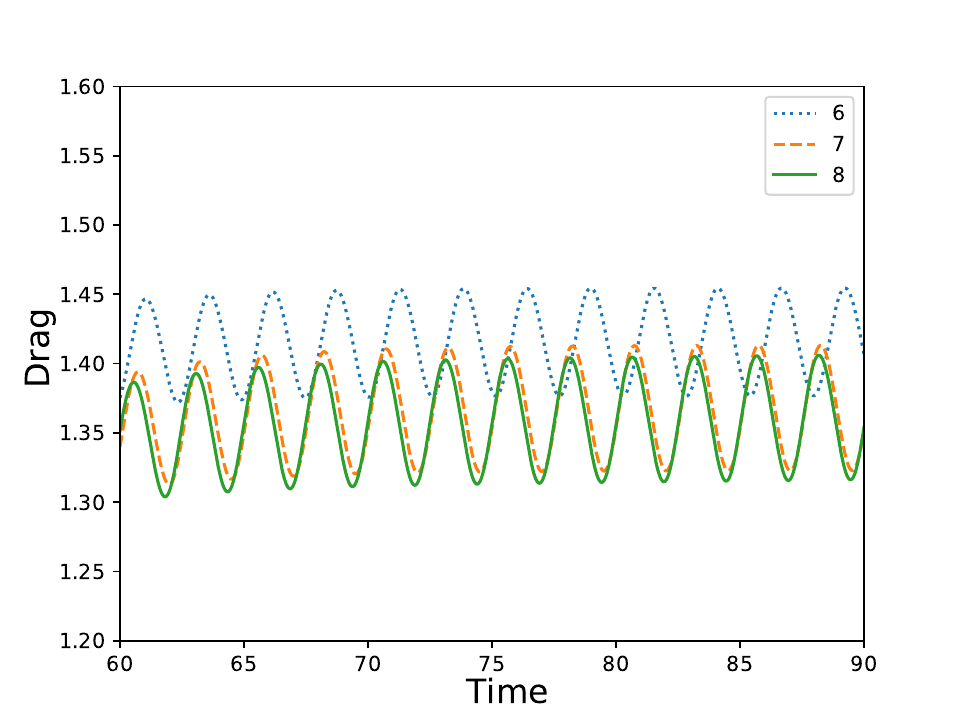}}
      \resizebox{200pt}{!}{\includegraphics{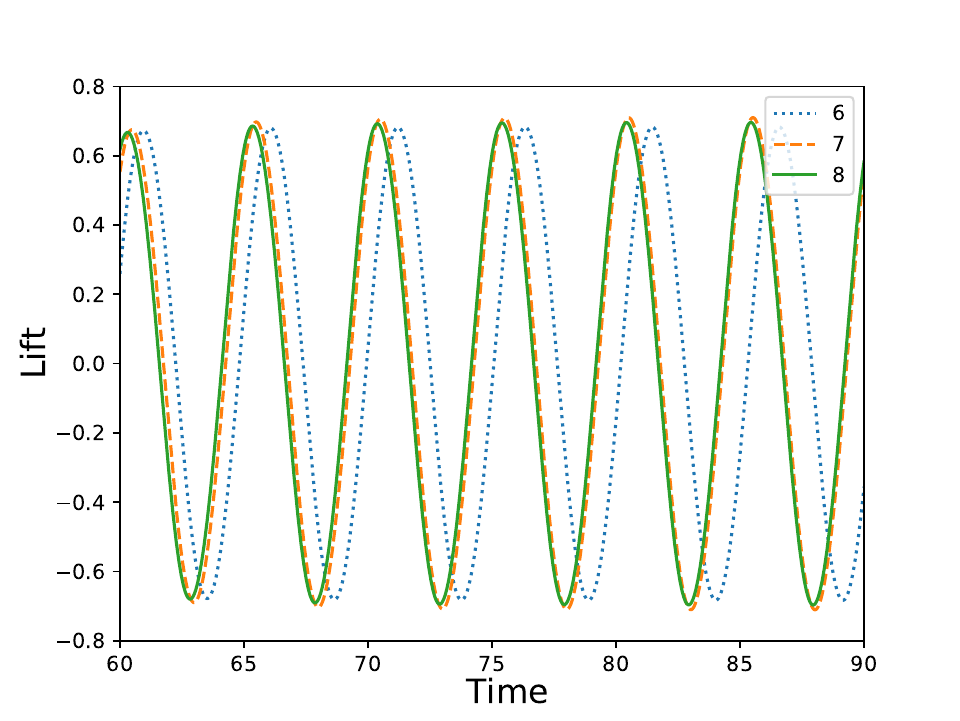}}
      
  \resizebox{200pt}{!}{\includegraphics{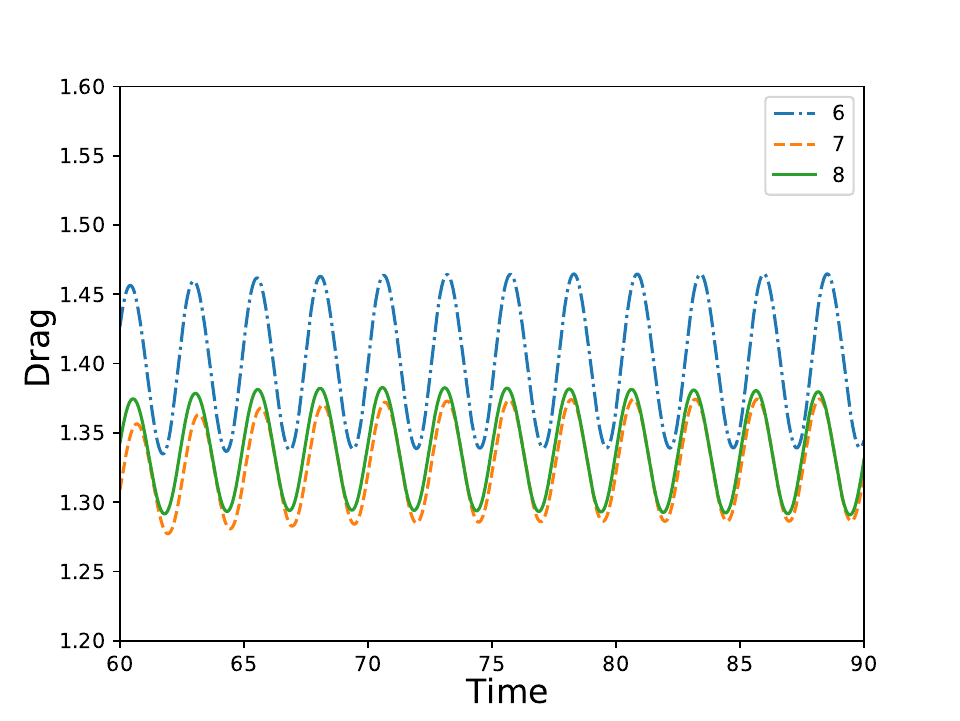}}
    \resizebox{200pt}{!}{\includegraphics{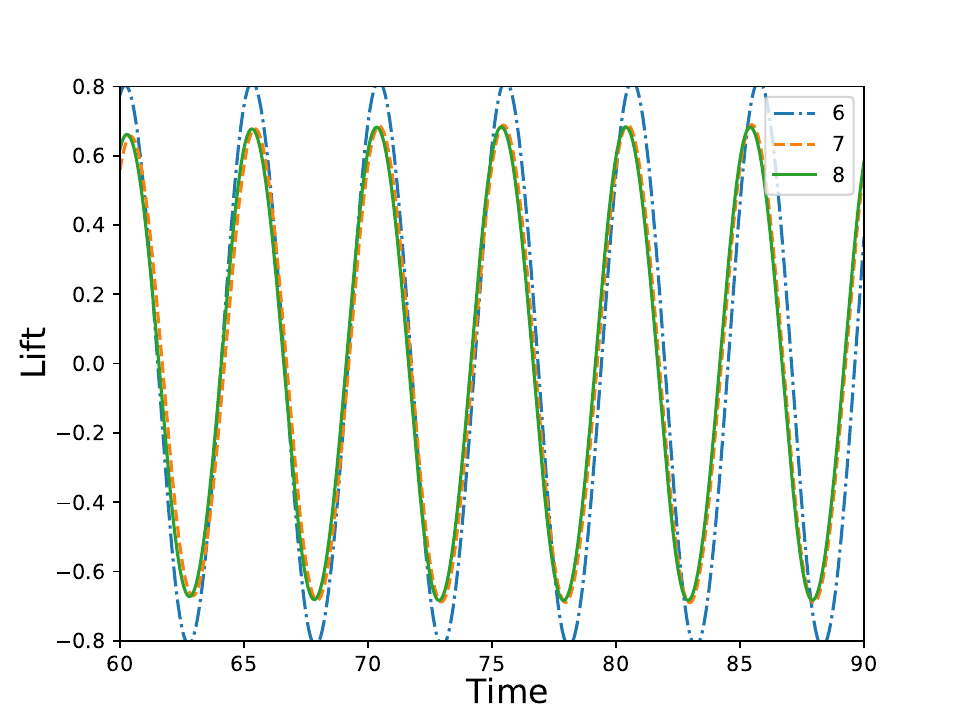}}
  \caption{Time histories of lift and drag for flow past a circular cylinder for $\ell_\text{max} = 6$, $7$, and $8$. From top to bottom, the normal vectors used are flat, inverse centroid-weighted, and projected. Both lift and drag exhibit convergence for each representation of the normal vector.}
\end{figure}

We quantify the simulation outputs using the nondimensionalized drag and lift coefficients, which are computed via
\begin{equation}
  C_{\text{D}} = \frac{- \int_{\Gamma_0} F^x (\mathbf{X}, t)
  \, {\mathrm d}A}{\frac{1}{2} \rho U^2 D}
\end{equation}
\begin{equation}
  {C_{\text{L}}}_{\textrm{}} = \frac{- \int_{\Gamma_0} F^y (\mathbf{X}, t)
  \, {\mathrm d}A}{\frac{1}{2} \rho U^2 D}
\end{equation}
As a verification of the method, we conduct a grid convergence test (Figure 4). For both lift and drag coefficients, all normal vector representations converge to quantities within the range of previous literature values. The smoothed representations of normal vectors converge more rapidly under grid refinement.
\iffalse
For a discrete
collection of times $\{ t_i \}_{i = 1}^M$ and the corresponding drag or lift
coefficients and $C (t_i)$, we compute the average coefficient as $C_{\text{avg}} = \frac{\max(C)+\min(C)}{2}$. Root-mean-squared drag and lift are computed as

\begin{equation}
C_{\text{rms}} = \sqrt{\frac{1}{M} \sum_{i = 1}^M \left(
C (t_i) - C_{\text{avg}} \right)^2}.
\end{equation}

\fi
The Strouhal number (St) is computed as $\text{St} = \frac{f D }{U}$, in which $f$ is the frequency of vortex shedding, $D$ is the diameter of the cylinder, and $U$ is the horizontal inflow velocity.

\begin{table}[H]
\center{}
  \begin{tabular}{l}
    \begin{tabular}{|c|c|c|c|}
      \hline
      \textbf{Reference} & $C_\text{D}$ & $C_\text{L}$ & St\\
      \hline
      Griffith and Luo{\cite{griffith_hybrid_2017}} & $1.360 \pm 0.046$ & $\pm
      0.70$ & 0.195\\
      \hline
      Liu et al. {\cite{liu_preconditioned_1998}} & $1.310 \pm 0.049$ & $\pm
      0.69$ & 0.192\\
      \hline
      Braza et al. {\cite{braza_numerical_1986}} & $1.400 \pm 0.050$ & $\pm
      0.75$ & 0.200\\
      \hline
      Calhoun{\cite{calhoun_cartesian_2002}} & $1.172 \pm 0.058$ & $\pm
      0.67$ & 0.202\\
      \hline
      Present Flat & $1.353 \pm 0.046$ &
      $\pm 0.70$ & 0.195\\
      \hline
      Present ICW & $1.360 \pm 0.047$ & $\pm 0.70$& 0.195\\
      \hline
      Present Project & $1.337 \pm 0.048$ & $\pm 0.69$& 0.197\\
      \hline
    \end{tabular} 
  \end{tabular}
  \caption{Comparisons of lift and drag coefficients for $\Rey = 200$ flow past a circular cylinder to results from prior studies.  }
\end{table}
The present lift and drag values fall within the
range of literature values (Table 1). Griffith and Luo\cite{griffith_hybrid_2017} use an immersed
boundary method with a finite element representation of the structure and a finite difference method for the incompressible Navier-Stokes equation. Liu et al.\cite{liu_preconditioned_1998}~use a
finite differences with a turbulence model. Braza et al.\cite{braza_numerical_1986} use finite
differences on the Navier Stokes equation specifically for cylinder flow.
Calhoun\cite{calhoun_cartesian_2002} uses an IIM with finite differences on a Cartesian grid.
\subsection{Pressure Loaded Cylinder in Two Spatial Dimensions}

This section considers a pressure loaded cylinder in two  spatial dimensions. The cylinder of diameter $D=1$ is centered at the origin the computational domain $\Omega=[-1,1]^2$. The fluid density and viscosity are $\rho=1$ and $\mu = 0.2$. On each boundary, we impose zero tangential traction and zero normal velocity. The computational domain is discretized with $N=16$, $32$, and $64$ grid cells in both directions for convergence studies. The time step size is proportional to the grid spacing, such that $\Delta t=\frac{h}{500}$. For this time step size and grid spacing, the CFL number is approximately $3\cdot10^{-4}$. To impose rigid body motion, we set $\kappa=\frac{3N}{2}\cdot 10^3 $ and $\eta=0$. Figure 5 shows a visualization of the interface and the pressure field.

\begin{figure}[H]
\centering
  \resizebox{250pt}{!}{\includegraphics{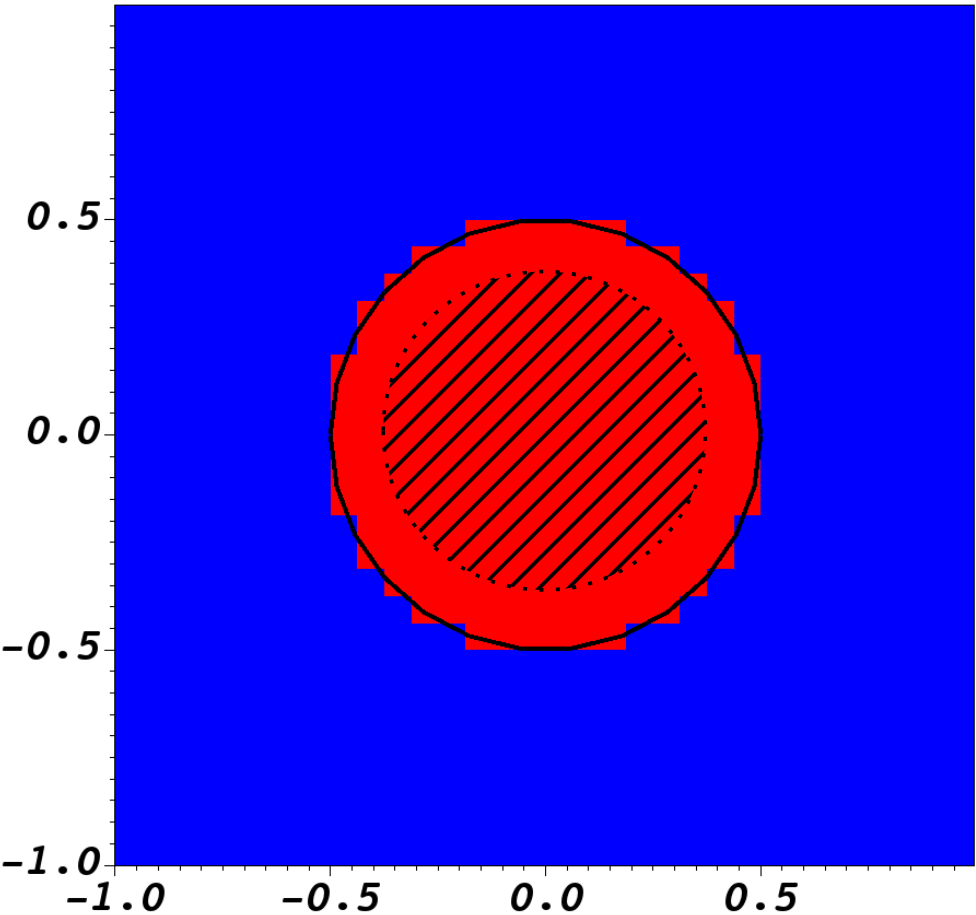}} 
  \resizebox{45pt}{235pt}{

\tikzset{every picture/.style={line width=0.75pt}} %set default line width to 0.75pt        

\begin{tikzpicture}[x=0.75pt,y=0.75pt,yscale=-1,xscale=1]
%uncomment if require: \path (0,384); %set diagram left start at 0, and has height of 384

%Image [id:dp3783569225831278] 
\draw (56.38,216.48) node  {\includegraphics[width=17.07pt,height=238.2pt]{colorbar_displacement.png}};
%Straight Lines [id:da6763773634376069] 
\draw [color={rgb, 255:red, 0; green, 0; blue, 0 }  ,draw opacity=1 ]   (67.76,58.68) -- (75.76,58.68) ;
%Straight Lines [id:da1207005768568743] 
\draw    (67.76,375.28) -- (76.76,375.28) ;
%Straight Lines [id:da48914525610221027] 
\draw    (67.94,215.76) -- (76.76,215.76) ;
%Straight Lines [id:da7189166377438487] 
\draw [color={rgb, 255:red, 0; green, 0; blue, 0 }  ,draw opacity=1 ]   (67.94,137.98) -- (76.76,137.98) ;
%Straight Lines [id:da7702940702509451] 
\draw    (67.94,292.54) -- (76.76,292.54) ;

% Text Node
\draw (49,33) node [anchor=north west][inner sep=0.75pt]  [font=\Large]  {$p$};
% Text Node
\draw (72.82,47.42) node [anchor=north west][inner sep=0.75pt]   [align=left] {100};
% Text Node
\draw (73.82,126.72) node [anchor=north west][inner sep=0.75pt]   [align=left] {75};
% Text Node
\draw (75.82,204.02) node [anchor=north west][inner sep=0.75pt]   [align=left] {50};
% Text Node
\draw (75.82,280.8) node [anchor=north west][inner sep=0.75pt]   [align=left] {25};
% Text Node
\draw (77.82,360.58) node [anchor=north west][inner sep=0.75pt]   [align=left] {0};

\end{tikzpicture}}
  \caption{The pressurized cylinder is visualized in the entire fluid domain. The shaded region inside the cylinder shows the region in which $Q$ is applied to impose $p=100$. We use $N=32$ grid cells in both directions.}
\end{figure}

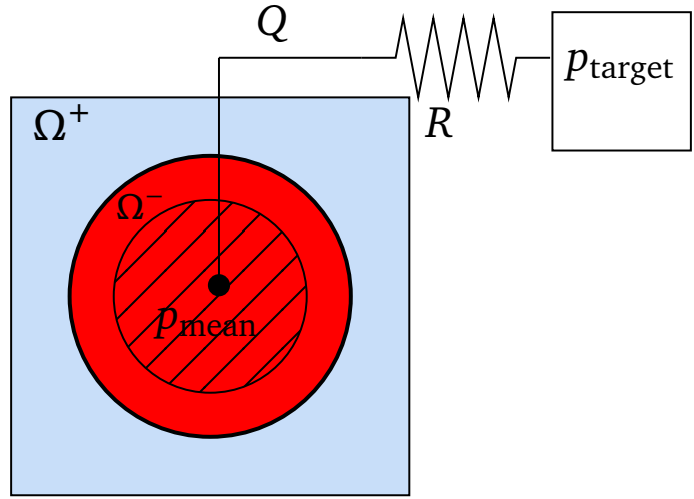
\begin{figure}[H]
\centering

% Pattern Info
 
\tikzset{
pattern size/.store in=\mcSize, 
pattern size = 5pt,
pattern thickness/.store in=\mcThickness, 
pattern thickness = 0.3pt,
pattern radius/.store in=\mcRadius, 
pattern radius = 1pt}
\makeatletter
\pgfutil@ifundefined{pgf@pattern@name@_f5nodo9jm}{
\pgfdeclarepatternformonly[\mcThickness,\mcSize]{_f5nodo9jm}
{\pgfqpoint{0pt}{0pt}}
{\pgfpoint{\mcSize+\mcThickness}{\mcSize+\mcThickness}}
{\pgfpoint{\mcSize}{\mcSize}}
{
\pgfsetcolor{\tikz@pattern@color}
\pgfsetlinewidth{\mcThickness}
\pgfpathmoveto{\pgfqpoint{0pt}{0pt}}
\pgfpathlineto{\pgfpoint{\mcSize+\mcThickness}{\mcSize+\mcThickness}}
\pgfusepath{stroke}
}}
\makeatother
\tikzset{every picture/.style={line width=0.75pt}} %set default line width to 0.75pt        

\begin{tikzpicture}[x=0.75pt,y=0.75pt,yscale=-1,xscale=1]
%uncomment if require: \path (0,300); %set diagram left start at 0, and has height of 300

%Shape: Square [id:dp7508677822825479] 
\draw  [fill={rgb, 255:red, 121; green, 176; blue, 242 }  ,fill opacity=0.41 ] (128,53) -- (328.82,53) -- (328.82,253.82) -- (128,253.82) -- cycle ;
%Shape: Circle [id:dp022756871865382666] 
\draw  [fill={rgb, 255:red, 255; green, 0; blue, 0 }  ,fill opacity=1 ][line width=1.5]  (157.5,153.41) .. controls (157.5,114.25) and (189.25,82.5) .. (228.41,82.5) .. controls (267.57,82.5) and (299.32,114.25) .. (299.32,153.41) .. controls (299.32,192.57) and (267.57,224.32) .. (228.41,224.32) .. controls (189.25,224.32) and (157.5,192.57) .. (157.5,153.41) -- cycle ;
%Shape: Circle [id:dp7991808988059292] 
\draw  [pattern=_f5nodo9jm,pattern size=14.774999999999999pt,pattern thickness=0.75pt,pattern radius=0pt, pattern color={rgb, 255:red, 0; green, 0; blue, 0}] (179.75,153.41) .. controls (179.75,126.54) and (201.54,104.75) .. (228.41,104.75) .. controls (255.28,104.75) and (277.07,126.54) .. (277.07,153.41) .. controls (277.07,180.28) and (255.28,202.07) .. (228.41,202.07) .. controls (201.54,202.07) and (179.75,180.28) .. (179.75,153.41) -- cycle ;
%Shape: Resistor [id:dp8974430244743337] 
\draw   (304.82,33) -- (321.92,33) -- (325.72,13) -- (333.32,53) -- (340.92,13) -- (348.52,53) -- (356.12,13) -- (363.72,53) -- (371.32,13) -- (378.92,53) -- (382.72,33) -- (399.82,33) ;
%Straight Lines [id:da32096932965398894] 
\draw    (304.82,33) -- (232.82,33) ;
%Straight Lines [id:da4336137463737789] 
\draw    (232.82,147.82) -- (232.82,33) ;
%Shape: Square [id:dp942021874118452] 
\draw   (401,10) -- (470.82,10) -- (470.82,79.82) -- (401,79.82) -- cycle ;
%Shape: Circle [id:dp10968830431381105] 
\draw  [fill={rgb, 255:red, 0; green, 0; blue, 0 }  ,fill opacity=1 ] (227.73,147.82) .. controls (227.73,145.01) and (230.01,142.73) .. (232.82,142.73) .. controls (235.63,142.73) and (237.91,145.01) .. (237.91,147.82) .. controls (237.91,150.63) and (235.63,152.91) .. (232.82,152.91) .. controls (230.01,152.91) and (227.73,150.63) .. (227.73,147.82) -- cycle ;

% Text Node
\draw (137,55.4) node [anchor=north west][inner sep=0.75pt]  [font=\LARGE]  {$\Omega ^{+}$};
% Text Node
\draw (179,100) node [anchor=north west][inner sep=0.75pt]  [font=\Large]  {$\Omega ^{-}$};
% Text Node
\draw (405,27.4) node [anchor=north west][inner sep=0.75pt]  [font=\LARGE]  {$p_{\text{target}}$};
% Text Node
\draw (198,157.4) node [anchor=north west][inner sep=0.75pt]  [font=\LARGE]  {$p_{\text{mean}}$};
% Text Node
\draw (250,5) node [anchor=north west][inner sep=0.75pt]  [font=\LARGE]  {$Q$};
% Text Node
\draw (335.32,56.4) node [anchor=north west][inner sep=0.75pt]  [font=\LARGE]  {$R$};

\end{tikzpicture}
\caption{A penalty source term $Q$ is applied within a sub-region of the cylinder to impose $p_\mathrm{target}$. The magnitude of $Q$ is determined by a feedback control loop analogous to a infinite reservoir at $p_{\mathrm{target}}$ connected via a hydraulic resistance.}
\label{fig:pressure_schematic}
\end{figure}

To establish the desired pressure $P_{\text{target}}$ within the cylinder, we introduce a volumetric source\cite{GRIFFITH200710,griffith2005simulating} term $Q(\mathbf{x},t)$ in the velocity divergence constraint, $\nabla \cdot \mathbf{u} = Q(\mathbf{x},t)$. 
This method is analogous to connecting the interior of the cylinder to an infinite reservoir at pressure $P_{\text{target}}$ through a connection with resistance $R$ (Figure~\ref{fig:pressure_schematic}).
The source term is radially smoothed using a cosine kernel, such that $Q(\mathbf{x},t) = Q(t) \cos\left(\frac{\pi r}{2 r_{\mathrm{src}}}\right)$ for $r < r_{\mathrm{src}}$ and zero otherwise.

\begin{figure}[H]
\centering
\resizebox{!}{200pt}{\includegraphics[]{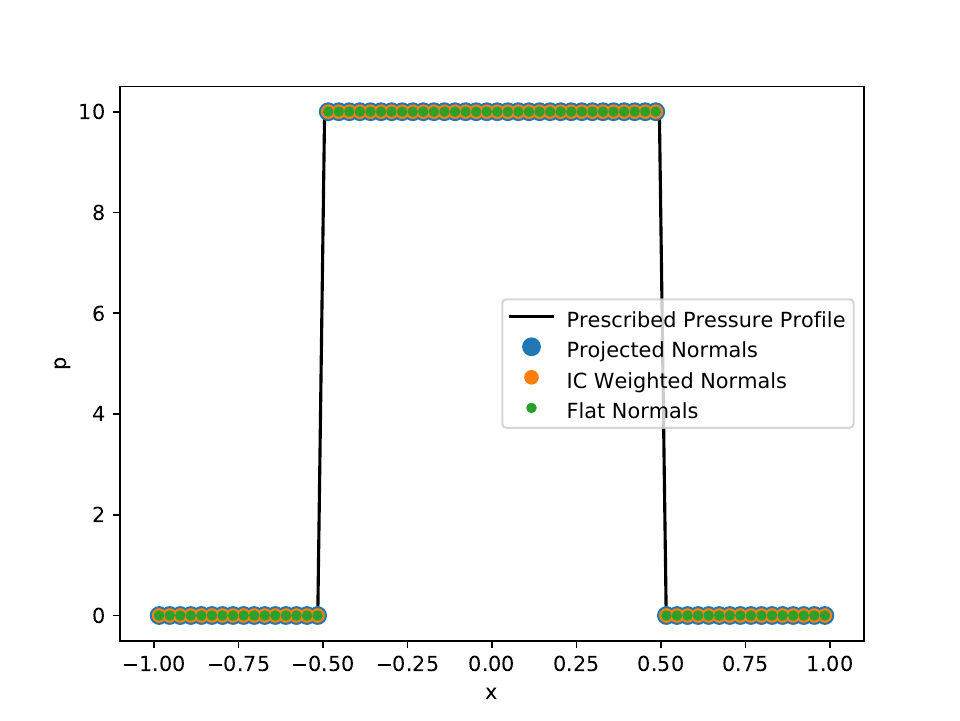}}

\caption{The pressure profile $p(x,0)$ of a pressurized cylinder in two spatial dimensions is visualized from $x=-1$ to $x=1.$ Smoothed normal vector representations of the jump condition use a PD controller, and flat normal representations use an ODE-based controller. Each of these converge to the imposed pressure profile. }
\label{fig:2d_pressure_profile}
\end{figure}

The source term is driven by a feedback control loop using two different methods. We use two methods for two reasons. The first method converges rapidly for smoothed normals, but flat normals exhibit enough leakage such that the method fails to converge to the target pressure. The second method provides a means to impose the target pressure for flat normals, but takes a burdensome amount of time to converge for smoothed normals. We use the first method for simulations using either inverse centroid-weighted or projected normals, and the second method for simulations using elemental normals. 
 %The choice of control strategy is tailored to the convergence properties of the normal vector representation.
 Using this match of surface normal representations and control strategies, we are able to rapidly converge $\tilde{p}(t)$ toward $p_\text{target}$. Figure~\ref{fig:2d_pressure_profile} demonstrates that each method used achieves the target pressure over the desired region.

The first method is a proportion-derivative (PD) controller based on the error between the target pressure and the instantaneous mean pressure $\bar{p}(t)$ measured within a source radius $r_{\mathrm{src}}$. 
The strength of the source term is updated at every time step according to:
\begin{equation}
Q(t) = \frac{1}{R} \left( e(t) +  \gamma \frac{\mathrm{d} e(t)}{\mathrm{d}t} \right).
\end{equation}
in which $e(t) = p_{\text{target}} - \bar{p}(t)$ is the pressure error, $\gamma=2\cdot10^{-3}$ is the derivative scaling parameter, and $R=450$ is the resistance. 
 The flow rate is adaptively adjusted until the internal pressure converges to $p_{\text{target}}$, at which point the system reaches a steady state where the net flow through the interface balances the source term.
 The second method models the pressure evolution as a first-order ordinary differential equation. In this formulation, the rate of change of the source term is proportional to the pressure mismatch,
 \begin{equation}
 \tau\frac{dQ(t)}{dt} = \frac{p_{\text{target}} - \bar{p}(t)}{R} .
 \end{equation}
 Discretized in time, the update rule becomes
 \begin{equation}
 Q^{n+1} = Q^n + \frac{\Delta t}{\tau R} e(t^n),
 \end{equation}
 in which $\tau = 5\Delta t$ is the time scaling parameter.

%  We utilize the feedback controller for centroid-weighted and consistent normals due to its rapid convergence, while the ODE method is employed for elemental normals to ensure efficient convergence.

\begin{figure}[H]
\centering
  \resizebox{250pt}{!}{\includegraphics{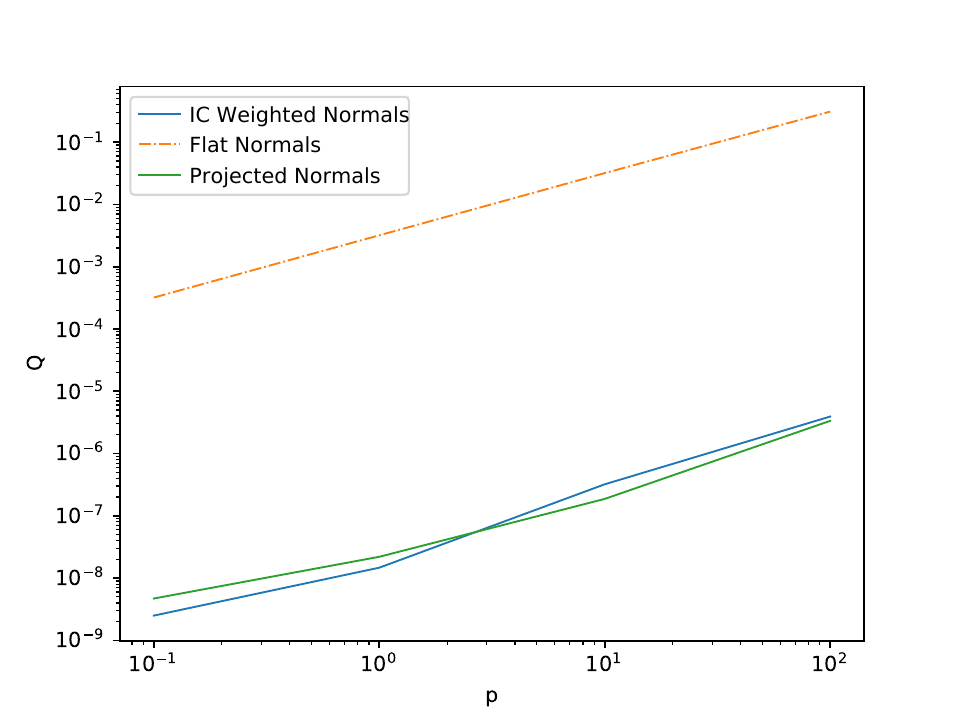}} 
  \caption{The steady state flow rates are visualized for various pressures $p_0=0.1$, $1$, $10$, and $100$ with $N=32$ grid cells in both directions. Using $L^2$ projected normals reduces numerical leakage by up to five orders of magnitude. }
  \label{fig:2d_various_pressures}
\end{figure}

To compare the accuracy of smoothed and flat normal vector representations, we pressurize the cylinder with various pressures $p_0 = 0.1$, $1$, $10$, and $100$ while fixing $N=32$ grid cells in both directions. 
We use penalty flow rate, $Q$, on a disk of radius 0.4 to impose the desired $p_0$ throughout the interior of the cylinder.
A steady state with $Q\neq0$ is interpretable as error in the interface's capacity to enforce the no-penetration condition. Figure~\ref{fig:2d_various_pressures} shows an accuracy improvement of about five orders of magnitude with smoothed interfacial normal vectors compared to elemental normal vectors.

\begin{figure}[H]
\centering
  \resizebox{250pt}{!}{\includegraphics{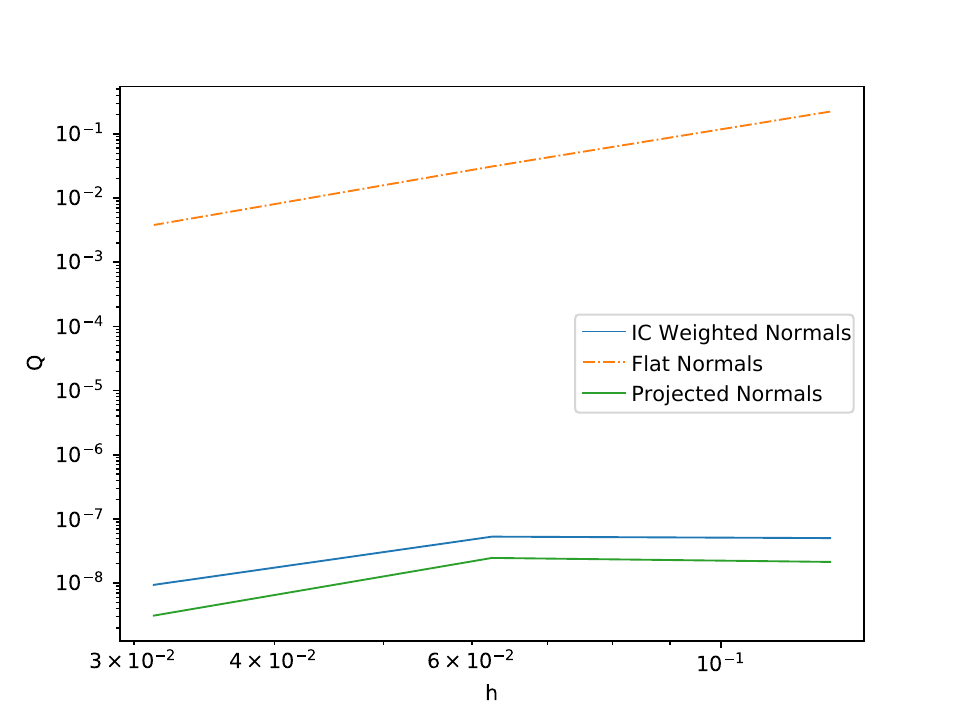}} 
  \caption{The steady state flow rates are visualized for various computational grid resolutions, $N=16$, $32$, and $64$ grid cells in both spatial directions. Using smoothed normals reduces numerical leakage by about six orders of magnitude at all grid resolutions. }
  \label{fig:2d_grid_conv}
\end{figure}

We also conduct a grid convergence study using a target pressure $p_0=10$ with $N=16$, $32$, and $64$ grid cells. Figure~\ref{fig:2d_grid_conv} shows the convergence of $Q$ under grid refinement. Using smoothed normal vectors results in a flow rate that is approximately $10^{-8}$ across all grid resolutions. Using smoothed normals improves volume conservation by at least to six orders of magnitude at all grid resolutions. 

\subsection{Pressure Loaded Cylinder in Three Spatial Dimensions}
\begin{figure}[H]
\centering
\hskip -.3cm
\includegraphics[width=.59\textwidth]{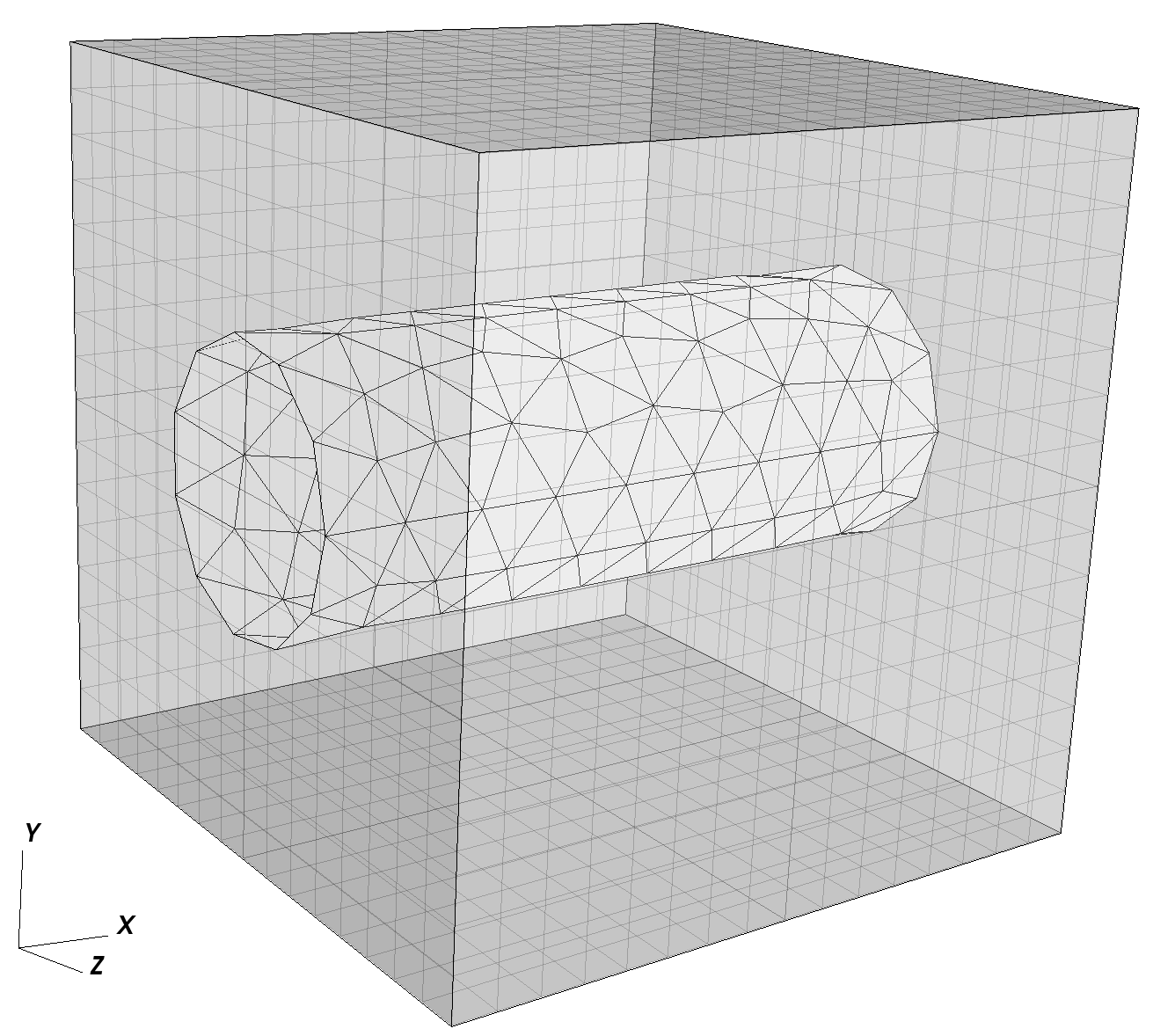}
  \caption{Computational meshes for flow in a horizontal cylinder, with $N = 16$ grid cells in each spatial direction. }
  %$h_{\Gamma}=0.5$,  $h=\frac{5}{16}$}
  \label{fig:mesh_tube}
\end{figure}

This section considers the interior flow through a rigid pressure-loaded cylinder of diameter $D=2$ and length $L=5$. The cylinder is parallel to the $x$-axis and centered at $y = z = 2.5$. 
The computational fluid domain is $\Omega = [0,5]^3$. Figure~\ref{fig:mesh_tube} shows the setup of the Eulerian and Lagrangian meshes in the fluid domain.
The fluid density and viscosity are $\rho = 1$ and $\mu = 0.04$.
On both the outflow and inflow boundary, we  impose $p = p_{0}$ inside the radius of the cylinder up to two mesh widths away from the interface.
Outside of this radius, we impose $u = v = w = 0$. On all other boundaries, we impose zero tangential traction and zero normal velocity.
The domain is discretized with a locally refined grid with $N=16$ cells on the coarsest level. The effective number of grid cells on the $\ell^{\text{th}}$ level is
${n_{\ell}}  = 2^{\ell - 1} n_{\text{c}}$. 
The time step size $\Delta t = 1.5\cdot 10^{-5}$.
For this time step size and grid spacing, the CFL number is approximately $6\cdot10^{-6}$. To keep the interface stationary, we set $\kappa = C_\kappa  N $ and $\eta = C_\eta N$, in which $C_\kappa = 15625$ and $C_\eta=312.5$.
To investigate the accuracy of the various normal vector representations, we measure the volumetric flow rate through the cylinder.
Imposing pressure loads of equal magnitude at either end of the cylinder should generate zero flow, allowing us to measure any flow as error.
\begin{figure}[H]
\centering
\hskip -.3cm
\includegraphics[width=.59\textwidth]{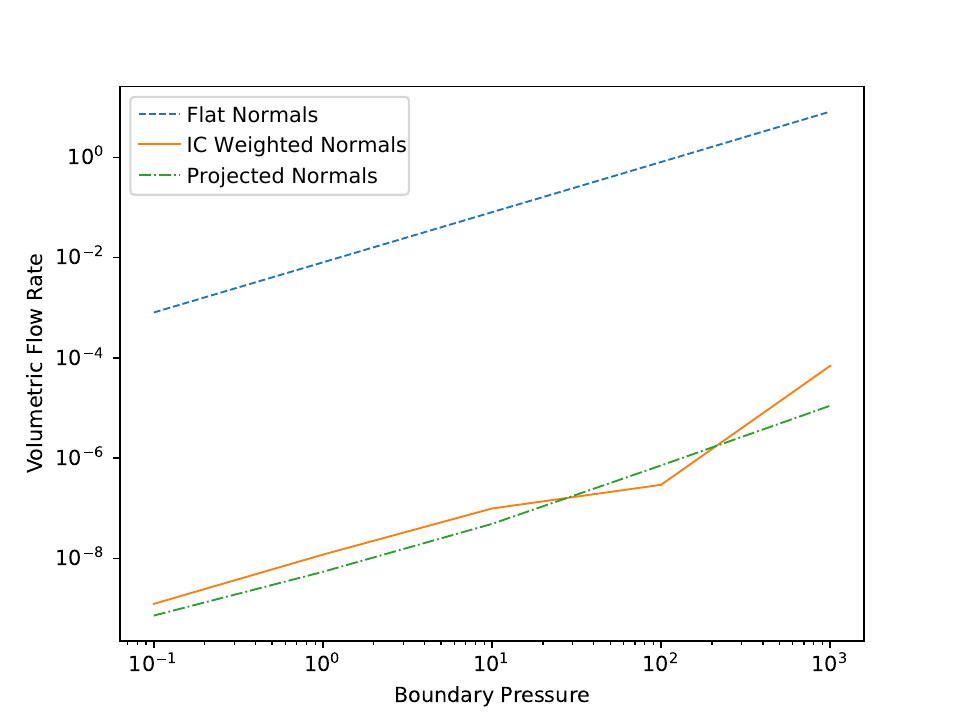}
  \caption{Comparison of the volumetric flow rates through a cylinder with boundary pressures $p_0 = 0.1$, $1$, $10$, $100$, and $1000$. }
  \label{fig:pressure_flowrates}
\end{figure}
Figure~\ref{fig:pressure_flowrates} shows the magnitude of the flow rate  compared to the boundary pressure $p_0 = 0.1$, $1$, $10$, $100$, and $1000$ for a fixed size of the Lagrangian mesh width $h_{\Gamma}=0.25$. At each pressure, the smoothed normal vector representation reduces the error by five to six orders of magnitude.

\begin{figure}[H]
\centering
\hskip -.3cm
\includegraphics[width=.59\textwidth]{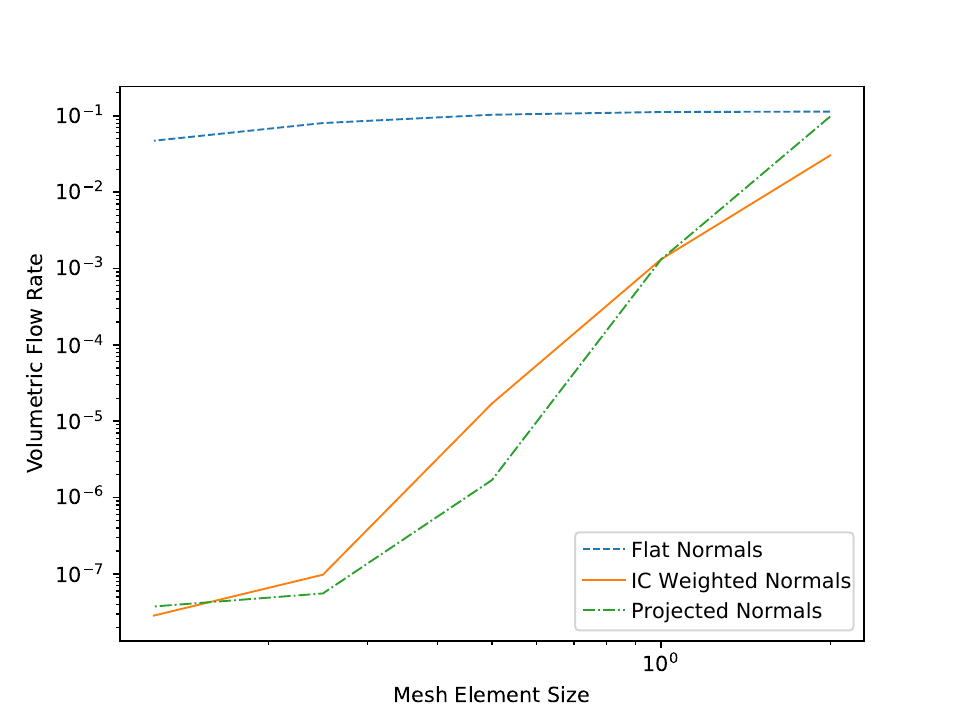}
  \caption{Comparison of the volumetric flow rates through a cylinder of Lagrangian mesh width $h_{\Gamma}=2$, $1$, $0.5$, $0.25$, and $0.125$ at a boundary pressure of $p_0 = 10$.}
  \label{fig:interface_sizes_flowrates}
\end{figure}
We also investigate the effects of Lagrangian mesh refinement on the volumetric flow rate. We fix the background grid and $p_0 = 10$, and we refine Lagrangian mesh such that $h_{\Gamma}=2$, $1$, $0.5$, and $0.25$ (Figure \ref{fig:interface_sizes_flowrates}). 
In the most under-resolved case in which the Lagrangian mesh is 12.8 times as coarse as the background grid, the three methods show similar error. 
In the most resolved case in which the coarsening ratio is 1.6, we observe approximately six orders of magnitude improvement while using the smoothed normal representations.

\begin{figure}[H]
\centering
\hskip -.3cm
\includegraphics[width=.59\textwidth]{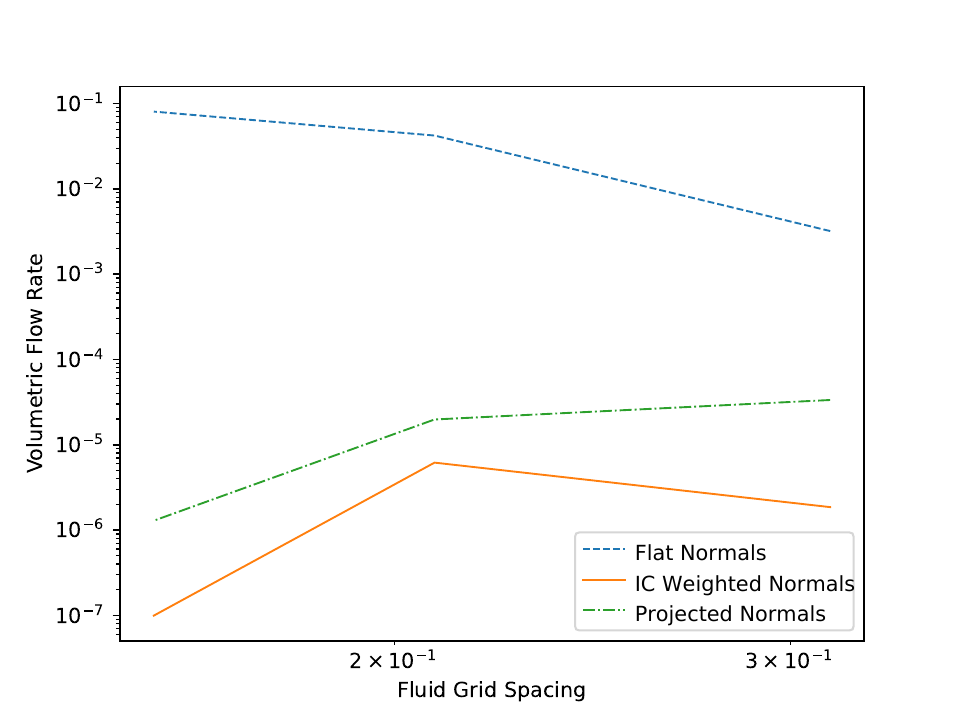}
  \caption{Comparison of the volumetric flow rates through a cylinder of Lagrangian mesh width $h_{\Gamma}=0.5$, $0.275$, and $0.25$ with fluid grid cell sizes $h = \frac{5}{16}$, $\frac{5}{24}$, and $\frac{5}{32}$ at a boundary pressure of $p_0 = 10$.}
  \label{fig:mesh_sizes_flowrates}
\end{figure}

We perform a grid convergence test by refining both the fluid grid and the interfacial mesh such that the coarsening ratio is fixed at 1.6. 
We use $N=16$, $24$, and $32$ fluid grid cells in each direction and interfacial mesh sizes $h_{\Gamma}=0.5$, $0.375$, and $0.25$ respectively (Figure~\ref{fig:mesh_sizes_flowrates}). 

\subsection{Poiseuille Flow in a Rigid Cylinder}
This section considers Poiseuille flow through a rigid pressure-loaded cylinder of diameter $D=2$ and length $L=5$. The cylinder is centered at $y=z=2.5$ and is parallel to the $x$-axis, extending from $x=0$ to $x=5$. 
The computational fluid domain is $\Omega = [0,5]^3$. The fluid density and viscosity are $\rho=1$ and $\mu=0.04$. 
On both the inflow and outflow boundaries, we impose $u=v=w=0$ outside the cylinder. 
Inside the cylinder, we set $p=p_{\text{inlet}}$ or $p=p_{\text{outlet}}$, and tangential velocities $v=w=0$. On all other boundaries, we impose zero normal velocity and zero tangential traction. The computational fluid domain is discretized into $N=32$ cells in each spatial direction. 

\begin{figure}[H]
\centering
  \resizebox{250pt}{!}{\includegraphics{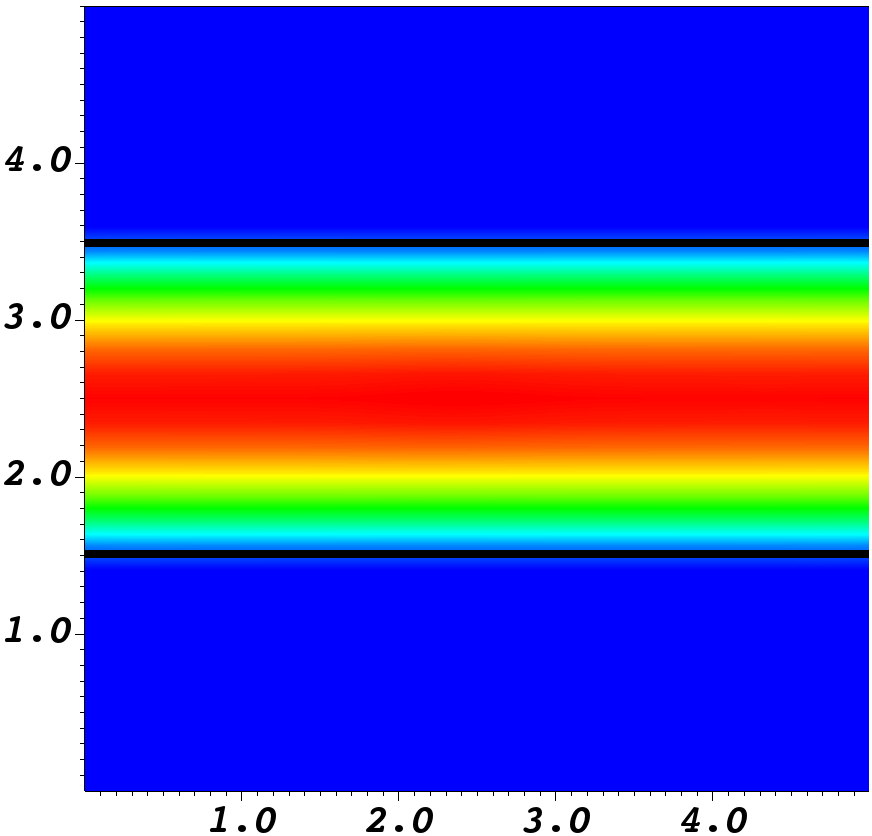}} 
  \resizebox{45pt}{200pt}{

\tikzset{every picture/.style={line width=0.75pt}} %set default line width to 0.75pt        

\begin{tikzpicture}[x=0.75pt,y=0.75pt,yscale=-1,xscale=1]
%uncomment if require: \path (0,384); %set diagram left start at 0, and has height of 384

%Image [id:dp3783569225831278] 
\draw (56.38,216.48) node  {\includegraphics[width=17.07pt,height=238.2pt]{colorbar_displacement.png}};
%Straight Lines [id:da6763773634376069] 
\draw [color={rgb, 255:red, 0; green, 0; blue, 0 }  ,draw opacity=1 ]   (67.76,58.68) -- (75.76,58.68) ;
%Straight Lines [id:da1207005768568743] 
\draw    (67.76,375.28) -- (76.76,375.28) ;
%Straight Lines [id:da48914525610221027] 
\draw    (67.94,215.76) -- (76.76,215.76) ;
%Straight Lines [id:da7189166377438487] 
\draw [color={rgb, 255:red, 0; green, 0; blue, 0 }  ,draw opacity=1 ]   (67.94,137.98) -- (76.76,137.98) ;
%Straight Lines [id:da7702940702509451] 
\draw    (67.94,292.54) -- (76.76,292.54) ;

% Text Node
\draw (49,33) node [anchor=north west][inner sep=0.75pt]  [font=\Large]  {$u$};
% Text Node
\draw (72.82,47.42) node [anchor=north west][inner sep=0.75pt]   [align=left] {1.262};
% Text Node
\draw (73.82,126.72) node [anchor=north west][inner sep=0.75pt]   [align=left] {0.946};
% Text Node
\draw (75.82,204.02) node [anchor=north west][inner sep=0.75pt]   [align=left] {0.631};
% Text Node
\draw (75.82,280.8) node [anchor=north west][inner sep=0.75pt]   [align=left] {0.315};
% Text Node
\draw (77.82,360.58) node [anchor=north west][inner sep=0.75pt]   [align=left] {8.2e-4};

\end{tikzpicture}}
  \caption{A slice of the $u$ velocity field is visualized over $\Omega =[0,5]^2$ in which $p_{\text{inlet}} = 1001$ and $p_{\text
{outlet}} = 1000$. A smooth parabolic profile is observed throughout the length of the cylinder.}
\label{fig:poiseuille_profile}
\end{figure}

We analyze two different pressure boundary condition configurations. In the first, we set $p_{\text{outlet}}=0$ and let $p_{\text{inlet}} = 0.1$, $1$, $10$, and $25$.
In this case, we set $\Delta t=7.5\cdot10^{-5}$. 
For this time step size and grid spacing, the CFL number is 0.005 once the model has reached steady state.
To impose the stationary boundary, we set $\kappa = 3\cdot10^3$ and $\eta = 5\cdot10^1$. In the second configuration, the axial pressure gradient is held constant while the pressure jump across the interface is varied, allowing for an isolated analysis of the radial pressure difference. To achieve this, we fix $p_\text{inlet} = p_\text{outlet} + 1$, and we vary $p_{\text{inlet}}= 1$, $10$, $100$, or $1000$.
In this case, we set $\Delta t =1\cdot10^{-6}$. For this time step size and grid spacing, the CFL number is approximately $8\cdot10^{-6}$. To impose rigid body motion, we set $\kappa = 2\cdot10^5$ and $\eta=1\cdot10^3$. 
Figure~\ref{fig:poiseuille_profile} shows a slice of the velocity field at steady state when $p_\text{inlet} = 1001$ and $p_\text{outlet} = 1000$ using smoothed normals. 

\begin{figure}[H]
\centering
\hskip -.3cm
\includegraphics[width=.59\textwidth]{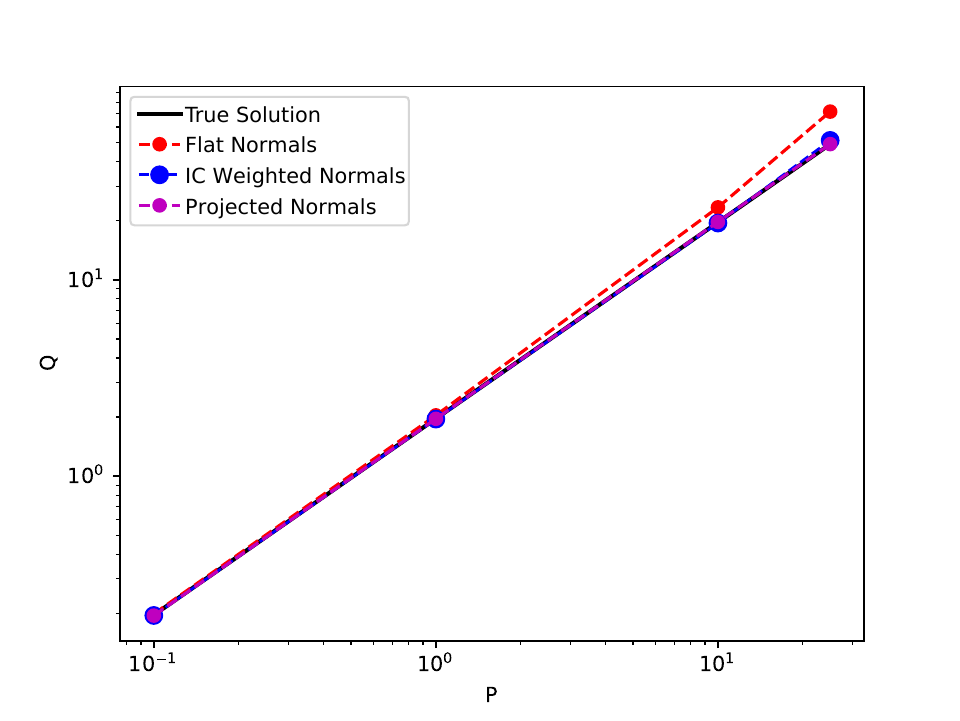}
  \caption{Comparison of the volumetric flow rates using both flat and smoothed normals in which $p_\text{inlet} = 0.1$, $1$, $10$, or $25$, and $p_\text{outlet} = 0$. Using smoothed normals produces a flow rate that agrees with the analytical rate, whereas a disparity arises while using flat normals at higher pressures.}
  \label{fig:poi_grad_diff}
\end{figure}

In both pressure configurations, we compare the numerical flow rate to the analytical flow rate using both smoothed and flat normals. The analytical flow rate is
\begin{equation}
Q = \frac{\pi R^4 \Delta p}{8 \mu L}.
\end{equation}
Figure~\ref{fig:poi_grad_diff} shows a comparison of the flow rates in the first pressure configuration.

We observe a greater error in flow rates at higher pressures when using flat normals. This error decreases when using smoothed normals. 

\begin{figure}[H]
\centering
\hskip -.3cm
\includegraphics[width=.59\textwidth]{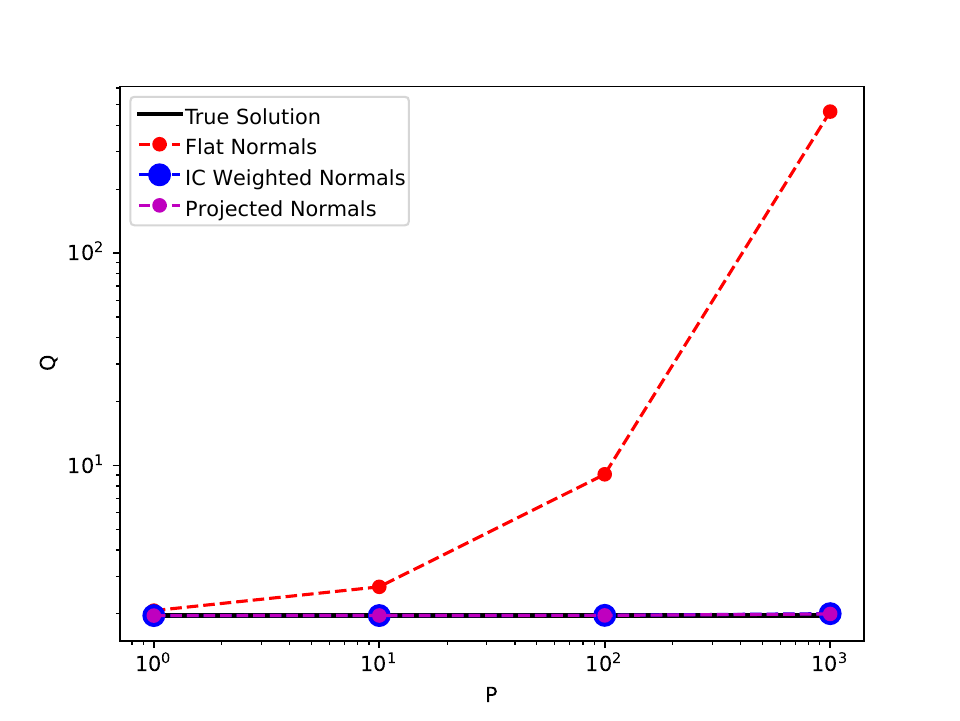}
  \caption{Comparison of the volumetric flow rates using both flat and smoothed normals in which $p_\text{inlet} = p_\text{outlet} + 1$, and $p_{\text{inlet}}= 1$, $10$, $100$, or $1000$. Using smoothed normals produces a flow rate that agrees with the analytical rate, whereas a disparity arises while using flat normals at higher pressures.}
  \label{fig:poi_grad_one}
\end{figure}

Figure~\ref{fig:poi_grad_one} shows a comparison of flow rates in the second pressure configuration.
We observe substantially larger error in the flow rate when using flat normals as the pressure difference across the interface increases. In contrast, the smoothed normal formulation agrees with the analytical flow rate across all tested pressures.

\section{Discussion}
In this work, we address a key limitation of the IIM, its failure to enforce the no-penetration condition across the interface, by introducing a new formulation for the projected interface condition based on discontinuous surface normals. Building on our earlier IIM framework, we present two practical reconstruction techniques. We use an $L^2$ projection and a computer-graphics inspired interpolation to construct globally continuous normal vector fields. 
We have demonstrated that capturing the curvature of the true geometry is crucial for accurately imposing jump conditions along the interface. 
These smoothing techniques substantially improve the no-penetration condition, by as much as six orders of magnitude, without sacrificing the method’s accuracy or robustness. 
Consequently, this advancement enables the IIM to surpass traditional immersed boundary-based approaches, which rely on  specific B-spline kernels or other treatments to mitigate the same leakage phenomena. Prior studies~\cite{GILMANOV2005457, MITTAL20084825 ,LEE20112677, SEO20117347} have noted the persistent problem of force oscillations in IB methods due to the pressure discontinuity across the immersed boundary, and have presented treatments to alleviate these numerical difficulties. These studies, however, have focused on bluff body flows as opposed to interior pressure loads.
The framework developed herein significantly broadens the applicability of the IIM, improving its ability to simulate complex, pressurized, three-dimensional flow problems in realistic settings.

\section*{Acknowledgments}
We gratefully acknowledge research support through National Institutes of Health awards R01HL157631 and R03HL182166
National Science Foundation award OAC 1931516. Computations were performed using
facilities provided by University of North Carolina at Chapel Hill through the Research Computing division of UNC Information Technology Services. We gratefully acknowledge Dr. Ebrahim M. Kolahdouz for providing foundational code that facilitated the setup of our numerical experiments. 

This manuscript is the result of funding in whole or in part by the National Institutes of Health (NIH). It is subject to the NIH Public Access Policy. Through acceptance of this federal funding, NIH has been given a right to make this manuscript publicly available in PubMed Central upon the Official Date of Publication, as defined by NIH.

\bibliographystyle{unsrt}
%\bibliography{literatur}
%\bibliography{biblio}

\end{document}